\documentclass[final]{siamltex}
\usepackage[letterpaper,text={6.5in,8.6in},centering]{geometry}
\usepackage{amssymb,amsmath,theorem}
\usepackage{subfigure,graphicx}

\newcommand{\norm}[1]{\ensuremath{\left\| #1 \right\|}}
\newcommand{\bracket}[1]{\ensuremath{\left[ #1 \right]}}

\newcommand{\parenth}[1]{\ensuremath{\left( #1 \right)}}
\newcommand{\refeqn}[1]{(\ref{eqn:#1})}
\newcommand{\reffig}[1]{Figure \ref{fig:#1}}

\newcommand{\deriv}[2]{\ensuremath{\frac{\partial #1}{\partial #2}}}
\newcommand{\SO}{\ensuremath{\mathrm{SO(3)}}}

\renewcommand{\Re}{\ensuremath{\mathbb{R}}}
\newcommand{\Sph}{\ensuremath{\mathbb{S}}}

\newtheorem{prop}{Proposition}
\newtheorem{cor}{Corollary}

{\theorembodyfont{\rmfamily}\newtheorem{example}{Example}}

\title{Lagrangian Mechanics and Variational Integrators on Two-Spheres\thanks{TL and ML have been supported in part by NSF Grant DMS-0504747 and DMS-0726263. TL and NHM have been supported in part by NSF Grant ECS-0244977 and CMS-0555797.}}

\author{Taeyoung Lee\thanks{Department of Aerospace Engineering, The University of Michigan, Ann Arbor, MI 48109 ({\tt tylee@umich.edu}).}%
    \and Melvin Leok\thanks{Department of Mathematics, Purdue University, West Lafayette, IN 47907 ({\tt mleok@math.purdue.edu}).}%
    \and N. Harris McClamroch\thanks{Department of Aerospace Engineering, The University of Michigan, Ann Arbor, MI 48109 ({\tt nhm@umich.edu}).}}

\begin{document}

\maketitle

\begin{abstract}
Euler-Lagrange equations and variational integrators are developed for Lagrangian mechanical systems evolving on a product of two-spheres. The geometric structure of a product of two-spheres is carefully considered in order to obtain global equations of motion. Both  continuous equations of motion and variational integrators completely avoid the singularities and complexities introduced by local parameterizations or explicit constraints. We derive global expressions for the Euler-Lagrange equations on two-spheres which are more compact than existing equations written in terms of angles. Since the variational integrators are derived from Hamilton's principle, they preserve the geometric features of the dynamics such as symplecticity, momentum maps, or total energy, as well as the structure of the configuration manifold. Computational properties of the variational integrators are illustrated for several mechanical systems.
\end{abstract}



\pagestyle{myheadings}
\thispagestyle{plain}
\markboth{T. Lee, M. Leok, and N. H. McClamroch}{Lagrangian Mechanics and Variational Integrators on Two-Spheres}

\section{Introduction}

The two-sphere $\Sph^2$ is the set of all points in the Euclidean space $\Re^3$ which are a unit distance from the origin. It is a two dimensional manifold that is locally diffeomorphic to $\Re^2$. Many classical and interesting mechanical systems, such as a spherical pendulum, a double spherical pendulum, and magnetic models, evolve on the two-sphere or on a product of two-spheres. In this paper, we derive Euler-Lagrange equations on configuration spaces of the form $(\Sph^2)^n$, for a positive integer $n$. We also develop geometric numerical integrators referred to as discrete Euler-Lagrange equations or variational integrators on $(\Sph^2)^n$.

In most of the literature that treats dynamic systems on $(\Sph^2)^n$, either $2n$ angles or $n$ explicit equality constraints enforcing unit length are used to describe the configuration of the system~\cite{BenSan.DDNS06,MarSchWen.ICIAM95}. These descriptions involve complicated trigonometric expressions and introduce additional complexity in analysis and computations. In this paper, we focus on developing continuous equations of motion and discrete equations of motion directly on $(\Sph^2)^n$, without need of local parameterizations, constraints, or reprojections. This provides a remarkably compact form of the equations of motion.

Geometric numerical integrators are numerical integration algorithms that preserve the geometric structure of the continuous dynamics, such as invariants, symplecticity, and the configuration manifold~\cite{Hai.BK00}. Conventional numerical integrators construct a discrete approximation to the flow using only information about the vector field, and ignore the physical laws and the geometric properties inherent in the differential equations~\cite{MclQui.BC01}. Consequently, they do not preserve important characteristics of the dynamics of the continuous equations of motion. In contrast, variational integrators are constructed by discretizing Hamilton's principle, rather than discretizing the continuous Euler-Lagrange equation~\cite{MosVes.CMP91,MaWe2001}. Since they are developed by using a discrete version of a physical principle, the resulting integrators have the desirable property that they are symplectic and momentum preserving, and they exhibit good energy behavior for exponentially long times.

Geometric numerical integration on $\Sph^2$ has been studied in \cite{MunZan.BC97,LewOlv.FCM03,LewNig.JCAM03}. The two-sphere is a homogeneous manifold; the special orthogonal group $\SO$ acts transitively on $\Sph^2$, and Lie group methods~\cite{IsMuNoZa2000} can be adapted to generate numerical flows on $\Sph^2$.

In this paper, we study Lagrangian mechanical systems on $(\Sph^2)^n$. Thus, it is desirable to preserve the geometric properties of the dynamics, such as momentum map, symplecticity, and total energy, in addition to the structure of the configuration manifold~\cite{CMDA07}. We combine the approaches of geometric integrators on homogeneous manifolds and variational integrators to obtain variational integrators on $(\Sph^2)^n$ that preserve the geometric properties of the dynamics as well as the homogeneous structure of the configuration manifold $(\Sph^2)^n$ concurrently.

The contributions of this paper can be summarized in two aspects. In the continuous time setting, the global Euler-Lagrange equations on $(\Sph^2)^n$ are developed in a compact form without local parameterization or constraints. This provides insight into the global dynamics on $(\Sph^2)^n$, which is desirable for theoretical studies. As a geometric numerical integrator, the discrete Euler-Lagrange equations on $(\Sph^2)^n$ are unique in the sense that they conserve both the geometric properties of the dynamics and the manifold structure of $(\Sph^2)^n$ simultaneously. The exact geometric properties of the discrete flow not only generate improved qualitative behavior, but they also provide accurate and reliable computational results in long-time simulation.

This paper is organized as follows. Lagrangian mechanics on $(\Sph^2)^n$ is described in Section \ref{sec:EL}. Variational integrators on $(\Sph^2)^n$ are developed in Section \ref{sec:DEL}. Computational properties are illustrated for several mechanical systems, namely a double spherical pendulum, an $n$-body problem on a sphere, an interconnected system of spherical pendula, pure bending of an elastic rod, a spatial array of magnetic dipoles, and  molecular dynamics that evolves on a sphere.

\section{Lagrangian mechanics on $(\Sph^2)^n$}\label{sec:EL}
In this section, continuous equations of motion for a mechanical system defined on $(\Sph^2)^n$ are developed in the context of Lagrangian mechanics. It is common in the published literature that the equations of motion are developed by using either two angles or a unit length constraint to characterize $\Sph^2$. Any description with two angles has singularities, and any trajectory near a singularity experiences numerical ill-conditioning. The unit length constraint leads to additional complexity in numerical computations. We develop global continuous equations of motion without resorting to local parameterizations or constraints. To achieve this, it is critical to understand the global characteristics of a mechanical system on $(\Sph^2)^n$. This section provides a good background for understanding the theory of discrete Lagrangian mechanics on $(\Sph^2)^n$ to be introduced in the next section.

The two-sphere is the set of points that have the unit length from the origin of $\Re^3$, i.e. $\Sph^2=\{q\in\Re^3\,|\, q\cdot q=1\}$. The tangent space $T_q\Sph^2$ for $q\in\Sph^2$ is a plane tangent to the two-sphere at the point $q$. Thus, a curve $q:\Re\rightarrow\Sph^2$ and its time derivative satisfy $q\cdot \dot q=0$.
The time-derivative of a curve can be written as
\begin{align}
    \dot q = \omega\times q,\label{eqn:dotq}
\end{align}
where the angular velocity $\omega\in\Re^3$ is constrained to be orthogonal to $q$, i.e. $q\cdot \omega=0$. The time derivative of the angular velocity is also orthogonal to $q$, i.e. $q\cdot\dot\omega=0$.

\subsection{Euler-Lagrange equations on $(\Sph^2)^n$}

We consider a mechanical system on the configuration manifold $\Sph^2\times\cdots \times \Sph^2=(\Sph^2)^n$. We assume that the Lagrangian $L:T(\Sph^2)^n\rightarrow\Re$ is given by the difference between a quadratic kinetic energy and a configuration-dependent potential energy as follows.
\begin{align}
    L(q_1,\ldots,q_n,\dot q_1,\ldots,\dot q_n) = \frac{1}{2} \sum_{i,j=1}^n  M_{ij}\dot q_i \cdot \dot q_j - V(q_1,\ldots,q_n),\label{eqn:L}
\end{align}
where $(q_i,\dot q_i)\in T\Sph^2$ for $i\in\{1,\ldots,n\}$, and $M_{ij}\in\Re$ is the $i,j$-th element of a symmetric positive definite inertia matrix $M\in\Re^{n\times n}$ for $i,j\in\{1,\ldots,n\}$. The configuration dependent potential is denoted by $V:(\Sph^2)^n\rightarrow\Re$.

The action integral is defined as the time integral of the Lagrangian, and the variation of the action integral leads to continuous equations of motion by applying Hamilton's principle. These are standard procedures to derive the Euler-Lagrange equations. The expression for the infinitesimal variation of $q_i\in\Sph^2$ should be carefully developed, since the configuration manifold is not a linear vector space. As in \refeqn{dotq}, the infinitesimal variation of $q_i$ can be written as a vector cross product,
\begin{align}
    \delta q_i = \xi_i \times q_i,\label{eqn:delqi}
\end{align}
where $\xi_i\in\Re^3$ is constrained to be orthogonal to $q_i$, i.e. $\xi_i\cdot q_i=0$. From this, the expression for the infinitesimal variation of $\dot q_i$ is given by
\begin{align}
    \delta \dot q_i = \dot\xi_i \times q_i + \xi_i\times \dot q_i.\label{eqn:delqidot}
\end{align}
These expressions are the key elements to obtaining global equations of motion on $(\Sph^2)^n$.

The variation of the Lagrangian can be written as
\begin{align*}
    \delta L = \sum_{i,j=1}^n  \delta \dot q_i \cdot M_{ij}\dot q_j - \sum_{i=1}^n \delta q_i \cdot \deriv{V}{q_i},
\end{align*}
where the symmetric property $M_{ij}=M_{ji}$ is used. Substituting \refeqn{delqi} and \refeqn{delqidot} into this, and using the vector identity $(a\times b)\cdot c = a\cdot (b\times c)$ for any $a,b,c\in\Re^3$, we obtain
\begin{align*}
    \delta L = \sum_{i,j=1}^n \dot \xi_i \cdot (q_i\times M_{ij}\dot q_j) + \xi_i\cdot(\dot q_i\times M_{ij}\dot q_j) - \sum_{i=1}^n \xi_i\cdot\parenth{q_i\times \deriv{V}{q_i}}.
\end{align*}
Let $\mathfrak{G}$ be the action integral defined as $\mathfrak{G}=\int_0^T L(q_1,\ldots,q_n,\dot q_1,\ldots,\dot q_n)\,dt$. Using the above equation and integrating by parts, the variation of the action integral is given by
\begin{align*}
    \delta \mathfrak{G} & = \sum_{i,j=1}^n \xi_i \cdot (\dot q_i\times M_{ij} \dot q_j + q_i\times M_{ij} \ddot q_j) \bigg|_{0}^T - \sum_{i=1}^n \int_0^T \xi_i \cdot \bracket{(q_i\times \sum_{j=1}^n M_{ij}\ddot q_j) + q_i\times \deriv{V}{q_i}}.
\end{align*}
From Hamilton's principle, $\delta \mathfrak{G}=0$ for any $\xi_i$ vanishing at $t=0,T$. Since $\xi_i$ is orthogonal to $q_i$, the continuous equations of motion satisfy
\begin{align}
    (q_i\times \sum_{j=1}^n M_{ij}\ddot q_j) + q_i\times \deriv{V}{q_i} = c_i(t) q_i\label{eqn:EL0}
\end{align}
for some scalar valued functions $c_i(t)$ for $i\in\{1,\ldots,n\}$. Taking the dot product of \refeqn{EL0} and $q_i$ implies that $0=c_i(t)\|q_i\|^2=c_i(t)$, which is to say that the scalar valued functions are uniformly zero. Now we find an expression for $\ddot q_i$. Since the left hand side expression is perpendicular to $q_i$, it is zero if and only if its cross product with $q_i$ is zero. Thus, we obtain
\begin{align}
    q_i\times (q_i\times \sum_{j=1}^n M_{ij}\ddot q_j) + q_i\times \parenth{q_i\times \deriv{V}{q_i}} = 0.\label{eqn:EL1}
\end{align}
From the vector identity $a\times (b\times c) = (a\cdot c)b- (a\cdot b) c$ for any $a,b,c\in\Re^3$, we have
\begin{align*}
    q_i\times (q_i\times \ddot q_i) & = (q_i\cdot\ddot q_i) q_i - (q_i\cdot q_i) \ddot q_i,\\
    & = - (\dot q_i \cdot \dot q_i)  q_i - \ddot q_i,
\end{align*}
where we use the properties $\frac{d}{dt}(q_i\cdot \dot q_i)= q_i\cdot\ddot q_i + \dot q_i\cdot \dot q_i=0$ and  $q_i\cdot q_i=1$. Substituting these into \refeqn{EL1}, we obtain an expression for $\ddot q_i$, which is summarized as follows.

\begin{prop}
Consider a mechanical system on $(\Sph^2)^n$ whose Lagrangian is expressed as \refeqn{L}. The continuous equations of motion are given by
\begin{align}
    M_{ii} \ddot q_i = q_i\times (q_i \times \sum_{\substack{j=1\\j\neq i}}^n M_{ij} \ddot q_j) -(\dot q_i\cdot \dot q_i)M_{ii}q_i + q_i \times \parenth{ q_i \times \deriv{V}{q_i}}\label{eqn:EL}
\end{align}
for $i\in\{1,\ldots,n\}$. Equivalently, this can be written in matrix form as
\begin{align}
    \begin{bmatrix}%
    M_{11}I_{3\times 3} & -M_{12} \hat q_1 \hat q_1 & \cdots & -M_{1n}\hat q_1 \hat q_1\\%
    -M_{21} \hat q_2\hat q_2 & M_{22} I_{3\times 3} & \cdots & -M_{2n} \hat q_2 \hat q_2\\%
    \vdots & \vdots & & \vdots\\
    -M_{n1} \hat q_n \hat q_n & -M_{n2}\hat q_n \hat q_n & \cdots & M_{nn} I_{3\times 3}
    \end{bmatrix}%
    \begin{bmatrix}
    \ddot q_1 \\ \ddot q_2 \\ \vdots \\ \ddot q_n
    \end{bmatrix}
    =
    \begin{bmatrix}
    -(\dot q_1 \cdot \dot q_1)M_{11} q_1 +\hat q_1^2 \deriv{V}{q_1}\\
    -(\dot q_2 \cdot \dot q_2)M_{22} q_2 +\hat q_2^2 \deriv{V}{q_2}\\
    \vdots\\
    -(\dot q_n \cdot \dot q_n)M_{nn} q_n +\hat q_n^2 \deriv{V}{q_n}
    \end{bmatrix},\label{eqn:ELm}
\end{align}
where the hat map $\hat\cdot : \Re^3\rightarrow \Re^{3\times 3}$ is defined such that $\hat a b = a\times b$ for any $a,b\in\Re^3$.
\end{prop}

Since $\dot q_i = \omega_i\times q_i$ for the angular velocity $\omega_i$ satisfying $q_i\cdot\omega_i=0$, we have
\begin{align*}
    \ddot q_i & = \dot \omega_i \times q_i + \omega_i\times (\omega_i\times q_i),\\
    & = \dot \omega_i \times q_i - (\omega_i\cdot\omega_i)q_i.
\end{align*}
Substituting this into \refeqn{EL0} and using the fact that $q_i\cdot \dot\omega_i=0$, we obtain continuous equations of motion in terms of the angular velocity.

\begin{cor}
The continuous equations of motion given by \refeqn{EL} can be written in terms of the angular velocity as
\begin{gather}
    M_{ii}\dot\omega_i = \sum_{\substack{j=1\\j\neq i}}^n \parenth{M_{ij} q_i\times(q_j\times \dot \omega_j) + M_{ij}(\omega_j\cdot\omega_j) q_i\times q_j} - q_i\times \deriv{V}{q_i},\label{eqn:ELw}\\
    \dot q_i = \omega_i\times q_i\label{eqn:ELq}
\end{gather}
for $i\in\{1,\ldots,n\}$. Equivalently, this can be written in matrix form as
\begin{align}
    \begin{bmatrix}%
    M_{11}I_{3\times 3} & -M_{12} \hat q_1 \hat q_2 & \cdots & -M_{1n}\hat q_1 \hat q_n\\%
    -M_{21} \hat q_2\hat q_1 & M_{22} I_{3\times 3} & \cdots & -M_{2n} \hat q_2 \hat q_n\\%
    \vdots & \vdots & & \vdots\\
    -M_{n1} \hat q_n \hat q_1 & -M_{n2}\hat q_n \hat q_2 & \cdots & M_{nn} I_{3\times 3}
    \end{bmatrix}%
    \begin{bmatrix}
    \dot \omega_1 \\ \dot \omega_2 \\ \vdots \\ \dot \omega_n
    \end{bmatrix}
    =
    \begin{bmatrix}
    \sum_{j=2}^n M_{1j}(\omega_j\cdot\omega_j)\hat q_1 q_j -\hat q_1\deriv{V}{q_1}\\
    \sum_{j=1,j\neq 2}^n M_{2j}(\omega_j\cdot\omega_j)\hat q_2 q_j -\hat q_2\deriv{V}{q_2}\\
    \vdots\\
    \sum_{j=1}^{n-1} M_{nj}(\omega_j\cdot\omega_j)\hat q_n q_j -\hat q_n\deriv{V}{q_n}\\
    \end{bmatrix}.\label{eqn:ELwm}
\end{align}
\end{cor}
Equations \refeqn{EL}--\refeqn{ELwm} are global continuous equations of motion for a mechanical system on $(\Sph^2)^n$. They avoid singularities completely, and they preserve the structure of $T(\Sph^2)^n$ automatically, if an initial condition is chosen properly. These equations are useful for understanding global characteristics of the dynamics. In addition, these expressions are dramatically more compact than the equations of motion written in terms of any local parameterization.

We need to check that the $3n\times 3n$ matrices given by the first terms of \refeqn{ELm} and \refeqn{ELwm} are nonsingular. This is a property of the mechanical system itself, rather than a consequence of the particular form of equations of motion. For example, when $n=2$, it can be shown that
\begin{align*}
    \det \begin{bmatrix} M_{11}I_{3\times 3} & -M_{12} \hat q_1 \hat q_1 \\
    -M_{12} \hat q_2 \hat q_2 & M_{22} I_{3\times 3} \end{bmatrix}
    & =  \det \begin{bmatrix} M_{11}I_{3\times 3} & -M_{12} \hat q_1 \hat q_2 \\
    -M_{12} \hat q_2 \hat q_1 & M_{22} I_{3\times 3} \end{bmatrix},\\
    & = M_{11}^2 M_{22}^2 ( M_{11}M_{22} - M_{12}^2 (q_1\cdot q_2)^2)( M_{11}M_{22} - M_{12}^2).
\end{align*}
Since the inertia matrix is symmetric positive definite, $M_{11},M_{22}>0$, $M_{11}M_{22} > M_{12}^2$, and from the Cauchy–-Schwarz inequality, $(q_1\cdot q_2)^2 \leq (q_1\cdot q_1)(q_2 \cdot q_2) =1$. Thus, the above matrices are non-singular. One may show a similar property for $n >2$. Throughtout this paper, it is assumed that the $3n\times 3n$ matrices given at the first terms of \refeqn{ELm} and \refeqn{ELwm} are nonsingular. Under this assumption, the Legendre transformation given in the next subsection is a diffeomorphism; the Lagrangian is hyperregular.

\subsection{Legendre transformation}

The Legendre transformation of the Lagrangian gives an equivalent Hamiltonian form of equations of motion in terms of conjugate momenta if the Lagrangian is hyperregular. Here, we find expressions for the conjugate momenta, which are used in the following section for the discrete equations of motion. For $q_i\in\Sph^2$, the corresponding conjugate momentum $p_i$ lies in the dual space $T_{q_i}^* \Sph^2$. We identify the tangent space $T_{q_i} \Sph^2$ and its dual space $T_{q_i}^* \Sph^2$ by using the usual dot product in $\Re^3$. The Legendre transformation is given by
\begin{align*}
    p_i \cdot \delta q_i & = \mathbf{D}_{\dot q_i} L(q_1,\ldots,q_n,\dot q_1,\ldots,\dot q_n)\cdot \delta q_i,\\
    & = \sum_{j=1}^n M_{ij} \dot q_j \cdot \delta q_i,
\end{align*}
which is satisfied for any $\delta q_i$ perpendicular to $q_i$. Here $\mathbf{D}_{\dot q_i}L$ denotes the derivative of the Lagrangian with respect to $\dot q_i$. The momentum $p_i$ is an element of the dual space identified with the tangent space, and the component parallel to $q_i$ has no effect since $\delta q_i\cdot q_i=0$. As such, the vector representing $p_i$ is perpendicular to $q_i$, and $p_i$ is equal to the projection of $\sum_{j=1}^n M_{ij} \dot{q}_j$ onto the orthogonal complement to $q_i$,
\begin{align}
    p_i & = \sum_{j=1}^n (M_{ij}\dot q_j - (q_i \cdot M_{ij}\dot q_j) q_i)= \sum_{j=1}^n ((q_i\cdot q_i)M_{ij}\dot q_j - (q_i \cdot M_{ij}\dot q_j) q_i),\nonumber\\
    & = M_{ii} \dot q_i -  q_i\times (q_i\times \sum_{\substack{j=1\\j\neq i}}^nM_{ij}\dot q_j).\label{eqn:pi}
\end{align}

\section{Variational integrators on $(\Sph^2)^n$}\label{sec:DEL}
The dynamics of Lagrangian and Hamiltonian systems on $(\Sph^2)^n$ have unique geometric properties; the Hamiltonian flow is symplectic, the total energy is conserved in the absence of non-conservative forces, and the momentum map associated with a symmetry of the system is preserved. The configuration space is a homogeneous manifold. These geometric features determine the qualitative dynamics of the system, and serve as a basis for theoretical study.

Conventional numerical integrators construct a discrete approximation of the flow using only information about the vector field. Other than the direction specified by the vector field, they completely ignore the physical laws and the geometric properties inherent in the differential equations~\cite{MclQui.BC01}. For example, if we integrate \refeqn{ELwm} by using an explicit Runge-Kutta method, the unit length of the vector $q_i$, and the total energy are not preserved numerically; we will see this later in this paper.

Numerical integration methods that preserve the simplecticity of a Hamiltonian system have been studied~\cite{San.AN92}. Coefficients of a Runge-Kutta method can be carefully chosen to satisfy a simplecticity criterion and order conditions to obtain a symplectic Runge-Kutta method. However, it can be difficult to construct such integrators, and it is not guaranteed that other invariants of the system, such as the momentum map, are preserved. Alternatively, variational integrators are constructed by discretizing Hamilton's principle, rather than by discretizing the continuous Euler-Lagrange equation~\cite{MosVes.CMP91,MaWe2001}. The key feature of variational integrators is that they are derived by a discrete version of a physical principle, so the resulting integrators satisfy the physical properties automatically in a discrete sense; they are symplectic and momentum preserving, and they exhibit good energy behavior for exponentially long times. Lie group methods are numerical integrators that preserve the Lie group structure of the configuration space~\cite{IsMuNoZa2000}. Recently, these two approaches have been unified to obtain Lie group variational integrators that preserve the geometric properties of the dynamics as well as the Lie group structure of the configuration manifold~\cite{CMA07}.

The two-sphere is a homogeneous manifold. It does not have a Lie group structure by itself, but instead, the special orthogonal group, $\SO=\{F\in\Re^{3\times 3}\,|\,F^T F=I_{3\times 3},\det F=1\}$, acts on $\Sph^2$ in a transitive way; for any $q_1,q_2\in\Sph^2$, there exist $F\in\SO$ such that $q_2=F q_1$. If a group acts transitively on a manifold, a curve on the manifold can be represented as the action of a curve in the Lie group on an initial point on the manifold. As such, Lie group methods can be applied to obtain numerical integration schemes for homogeneous manifolds~\cite{MunZan.BC97,LewNig.JCAM03,LewOlv.FCM03}. However, it is not guaranteed that these methods preserve the geometric properties of the dynamics. In this paper, we focus on a Lagrangian mechanical system evolving on the homogeneous manifold, $(\Sph^2)^n$ by extending the method of Lie group variational integrators~\cite{CMDA07,CMA07}. The resulting integrator preserves the dynamic characteristics and the homogeneous manifold structure concurrently.

\subsection{Discrete Euler-Lagrange equations on $(\Sph^2)^n$}

The procedure to derive discrete Euler-Lagrange equations follows the development of the continuous time case; the tangent bundle is replaced by a cartesian product of the configuration manifold, a discrete Lagrangian is chosen to approximate the integral of the Lagrangian over a discrete time step, and the variation of the corresponding discrete action sum provides discrete Euler-Lagrange equations, referred to as a variational integrator. The discrete version of the Legendre transformation yields the discrete equations in Hamiltonian form.

Let the number of timesteps be $N$, with constant timesteps $h>0$. A variable with subscript $k$ denotes the value of variable at $t=kh$. Define a discrete Lagrangian $L_d:(\Sph^2)^n\times (\Sph^2)^n\rightarrow\Re$ such that it approximates the integral of the Lagrangian given by \refeqn{L} over a discrete time step
\begin{align}
    L_d(q_{1_k},\ldots,q_{n_k},q_{1_{k+1}},\ldots,q_{n_{k+1}})%
     & = \frac{1}{2h} \sum_{i,j=1}^n M_{ij}(q_{i_{k+1}}-q_{i_k})\cdot(q_{j_{k+1}}-q_{j_k})%
     -\frac{h}{2} V_k - \frac{h}{2} V_{k+1},
\end{align}
where $V_k$ denotes the value of the potential at the $k$-th step, i.e. $V_k=V(q_{1_k},\ldots,q_{n_k})$. As given in \refeqn{delqi}, the infinitesimal variation of $q_{i_k}$ is written as
\begin{align}
    \delta q_{i_k} = \xi_{i_k} \times q_{i_k},\label{eqn:delqik}
\end{align}
where $\xi_{i_k}\in\Re^3$ is constrained to be orthogonal to $q_{i_k}$, i.e. $\xi_{i_k}\cdot q_{i_k}=0$. The variation of the discrete Lagrangian can be written as
\begin{align}
    \delta L_{d_k} = \frac{1}{h} \sum_{i,j=1}^n (\delta q_{i_{k+1}}-\delta q_{i_k} )\cdot M_{ij}(q_{j_{k+1}}-q_{j_k})-\frac{h}{2}\sum_{i=1}^n \parenth{\delta q_{i_k}\cdot \deriv{V_k}{q_{i_k}}+\delta q_{i_{k+1}}\cdot \deriv{V_{k+1}}{q_{i_{k+1}}}}.\label{eqn:Ldk0}
\end{align}
Substituting \refeqn{delqik} into \refeqn{Ldk0}, and using the vector identity $(a\times b)\cdot c = a\cdot (b\times c)$ for any $a,b,c\in\Re^3$, we obtain
\begin{align}
    \delta L_{d_k}  & = \frac{1}{h}\sum_{i,j=1}^n%
     \parenth{\xi_{i_{k+1}} \cdot (q_{i_{k+1}}\times M_{ij}(q_{j_{k+1}}-q_{j_k}))%
     -\xi_{i_{k}} \cdot (q_{i_{k}}\times     M_{ij}(q_{j_{k+1}}-q_{j_k}))}\nonumber\\%
     &\quad -\frac{h}{2}\sum_{i=1}^n \parenth{\xi_{i_k}\cdot\parenth{q_{i_k}\times \deriv{V_k}{q_{i_k}}}%
     +\xi_{i_{k+1}}\cdot\parenth{q_{i_{k+1}}\times \deriv{V_{k+1}}{q_{i_{k+1}}}}}.\label{eqn:Ldk}
\end{align}

Let $\mathfrak{G}_d$ be the discrete action sum defined as $\mathfrak{G}_d=\sum_{k=0}^{N-1} L_{d_k}$, which approximates the action integral as the discrete Lagrangian approximates a piece of the action integral over a discrete time step. The variation of the action sum is obtained by using \refeqn{Ldk}. Using the fact that $\xi_{i_k}$ vanish at $k=0$ and $k=N$, we can reindex the summation, which is the discrete analog of integration by parts, to yield
\begin{align*}
    \delta \mathfrak{G}_d =\sum_{k=1}^{N-1}\sum_{i=1}^n \xi_{i_{k}} \cdot\bracket{  %
    \frac{1}{h}(q_{i_{k}}\times \sum_{j=1}^n M_{ij}(-q_{j_{k+1}}+2q_{j_k}-q_{j_{k-1}}))%
    -hq_{i_k}\times \deriv{V_k}{q_{i_k}}}.
\end{align*}
From discrete Hamilton's principle $\delta\mathfrak{G}_d=0$ for any $\xi_{i_k}$ perpendicular to $q_{i_k}$. Using the same argument given in \refeqn{EL0}, the discrete equations of motion are given by
\begin{align}
    \frac{1}{h}(q_{i_{k}}\times \sum_{j=1}^n M_{ij}(-q_{j_{k+1}}+2q_{j_k}-q_{j_{k-1}}))%
    -hq_{i_k}\times \deriv{V_k}{q_{i_k}}=0\label{eqn:DEL0}
\end{align}
for $i\in\{1,\ldots\,n\}$. In addition, we require that the unit length of the vector $q_{i_k}$ is preserved. This is achieved by viewing $\Sph^2$ as a homogeneous manifold. Since the special orthogonal group $\SO$ acts on $\Sph^2$ transitively, we can define a discrete update map for $q_{i_k}$ as
\begin{align*}
    q_{i_{k+1}}= F_{i_k} q_{i_k}
\end{align*}
for $F_{i_k}\in\SO$. Then, the unit length of the vector $q_i$ is preserved through the discrete equations of motion, since  $q_{i_{k+1}}\cdot q_{i_{k+1}}=q_{i_k}^T F_{i_k}^T F_{i_k} q_{i_k}=1$. These results are summarized as follows.

\begin{prop}
Consider a mechanical system on $(\Sph^2)^n$ whose Lagrangian is expressed as \refeqn{L}. The discrete equations of motion are given by
\begin{gather}
    M_{ii} q_{i_{k}}\times F_{i_k} q_{i_k}
    +q_{i_{k}}\times \sum_{\substack{j=1\\j\neq i}}^n M_{ij} (F_{j_{k}}-I_{3\times 3}) q_{j_k}
    =
    q_{i_k}\times \sum_{j=1}^n M_{ij}(q_{j_k}-q_{j_{k-1}})%
    -h^2 q_{i_k}\times \deriv{V_k}{q_{i_k}}=0,\label{eqn:DEL}\\
    q_{i_{k+1}}= F_{i_k} q_{i_k}\label{eqn:qikp}
\end{gather}
for $i\in\{1,\ldots\,n\}$. For given $(q_{i_{k-1}},q_{i_k})$, we solve \refeqn{DEL} to obtain $F_{i_k}\in\SO$. Then, $q_{i_{k+1}}$ is computed by \refeqn{qikp}. This yields a discrete flow map $(q_{i_{k-1}},q_{i_k})\mapsto(q_{i_{k}},q_{i_{k+1}})$, and this process is repeated. 
\end{prop}

\subsection{Discrete Legendre transformation}
We find discrete equations of motion in terms of the angular velocity. The discrete Legendre transformation is given as follows~\cite{MaWe2001}.
\begin{align*}
    p_{i_k}\cdot \delta q_{i_k} & = -\mathbf{D}_{q_{i_k}} L_{d_k}\cdot \delta q_{i_k},\\
    & = \bracket{\frac{1}{h} \sum_{j=1}^n M_{ij} (q_{j_{k+1}}-q_{j_k})  + \frac{h}{2} \deriv{V_k}{q_{i_k}}}\cdot \delta q_{i_k},
\end{align*}
which can be directly obtained from \refeqn{Ldk0}. This is satisfied for any $\delta q_{i_k}$ perpendicular to $q_{i_k}$. Using the same argument used to derive \refeqn{pi}, the conjugate momenta $p_{i_k}$ is the projection of the expression in brackets onto the orthogonal complement of $q_{i_k}$. Thus, we obtain
\begin{align*}
    q_{i_k} = -\frac{1}{h} q_{i_k}\times (q_{i_k}\times \sum_{j=1}^n M_{ij} (q_{j_{k+1}}-q_{j_k})) - \frac{h}{2}q_{i_k}\times \parenth{q_{i_k}\times \deriv{V_k}{q_{i_k}}}.
\end{align*}
Comparing this to \refeqn{pi}, substituting $\dot q_{i_k}=\omega_{i_k}\times q_{i_k}$, and rearranging, we obtain
\begin{align*}
    q_{j_k}\times \bracket{M_{ii} \omega_{i_k} +   (q_{i_k}\times \sum_{\substack{j=1\\j\neq i}}^nM_{ij}(\omega_{j_k}\times q_{j_k}))-\frac{1}{h} (q_{i_k}\times \sum_{j=1}^n M_{ij} (q_{j_{k+1}}-q_{j_k})) - \frac{h}{2}q_{i_k}\times \deriv{V_k}{q_{i_k}}}=0.
\end{align*}
Since the expression in the brackets is orthogonal to $q_{i_k}$, the left hand side is equal to zero if and only if the expression in the brackets is zero. Thus,
\begin{align}
    M_{ii} \omega_{i_k} +   (q_{i_k}\times \sum_{\substack{j=1\\j\neq i}}^nM_{ij}(\omega_{j_k}\times q_{j_k}))=\frac{1}{h} (q_{i_k}\times \sum_{j=1}^n M_{ij} (q_{j_{k+1}}-q_{j_k})) + \frac{h}{2}q_{i_k}\times \deriv{V_k}{q_{i_k}}.\label{eqn:LT-}
\end{align}
This provides a relationship between $(q_{i_k},\omega_{i_k})$ and $(q_{i_k},q_{i_{k+1}})$. Comparing this with \refeqn{DEL0}, we obtain
\begin{align}
    M_{ii} \omega_{i_k} +   (q_{i_k}\times \sum_{\substack{j=1\\j\neq i}}^nM_{ij}(\omega_{j_k}\times q_{j_k}))=\frac{1}{h} (q_{i_k}\times \sum_{j=1}^n M_{ij} (q_{j_{k}}-q_{j_{k-1}})) - \frac{h}{2}q_{i_k}\times \deriv{V_k}{q_{i_k}},
    \label{eqn:LT+}
\end{align}
which provides a relationship between $(q_{i_k},\omega_{i_k})$ and $(q_{i_{k-1}},q_{i_{k}})$. Equations \refeqn{LT-} and \refeqn{LT+} give a discrete flow map in terms of the angular velocity; for a given $(q_{i_k},\omega_{i_k})$, we find $(q_{i_k},q_{i_{k+1}})$ by using \refeqn{LT-}. Substituting this into \refeqn{LT+} expressed at the $k+1$th step, we obtain $(q_{i_{k+1}},\omega_{i_{k+1}})$. This procedure is summarized as follows.

\begin{cor}
The discrete equations of motion given by \refeqn{DEL} and \refeqn{qikp} can be written in terms of the angular velocity as
\begin{gather}
    M_{ii} q_{i_{k}}\times F_{i_k} q_{i_k}%
    +q_{i_{k}}\times \sum_{\substack{j=1\\j\neq i}}^n M_{ij} (F_{j_{k}}-I_{3\times 3}) q_{j_k}%
    =M_{ii} h\omega_{i_k} -   (q_{i_k}\times \sum_{\substack{j=1\\j\neq i}}^nM_{ij}( q_{j_k}\times h\omega_{j_k}))%
    - \frac{h^2}{2}q_{i_k}\times \deriv{V_k}{q_{i_k}},\label{eqn:DELw1}\\
    q_{i_{k+1}}= F_{i_k} q_{i_k}\label{eqn:qikpw},\\
    M_{ii} \omega_{i_{k+1}} -   (q_{i_{k+1}}\times \sum_{\substack{j=1\\j\neq i}}^nM_{ij}( q_{j_{k+1}}\times \omega_{j_{k+1}}))=\frac{1}{h} (q_{i_{k+1}}\times \sum_{j=1}^n M_{ij} (q_{j_{k+1}}-q_{j_{k}})) - \frac{h}{2}q_{i_{k+1}}\times \deriv{V_{k+1}}{q_{i_{k+1}}}\label{eqn:DELw2}
\end{gather}
for $i\in\{1,\ldots,n\}$. Equivalently, \refeqn{DELw2} can be written in a matrix form as
\begin{align}
    &\begin{bmatrix}%
    M_{11}I_{3\times 3} & -M_{12} \hat q_{1_{k+1}} \hat q_{2_{k+1}} & \cdots & -M_{1n}\hat q_1 \hat q_{n_{k+1}}\\%
    -M_{21} \hat q_{2_{k+1}}\hat q_{1_{k+1}} & M_{22} I_{3\times 3} & \cdots & -M_{2n} \hat q_{2_{k+1}} \hat q_{n_{k+1}}\\%
    \vdots & \vdots & & \vdots\\
    -M_{n1} \hat q_{n_{k+1}} \hat q_{1_{k+1}} & -M_{n2}\hat q_{n_{k+1}} \hat q_{2_{k+1}} & \cdots & M_{nn} I_{3\times 3}
    \end{bmatrix}%
    \begin{bmatrix}
    \omega_{1_{k+1}} \\ \omega_{2_{k+1}} \\ \vdots \\ \omega_{n_{k+1}}
    \end{bmatrix}\nonumber\\
    &\qquad\qquad\qquad\qquad\qquad\qquad=\begin{bmatrix}
    \frac{1}{h} (q_{1_{k+1}}\times \sum_{j=1}^n M_{1j} (q_{j_{k+1}}-q_{j_{k}})) - \frac{h}{2}q_{1_{k+1}}\times \deriv{V_{k+1}}{q_{1_{k+1}}}\\
    \frac{1}{h} (q_{2_{k+1}}\times \sum_{j=1}^n M_{2j} (q_{j_{k+1}}-q_{j_{k}})) - \frac{h}{2}q_{2_{k+1}}\times \deriv{V_{k+1}}{q_{2_{k+1}}}\\
    \vdots\\
    \frac{1}{h} (q_{n_{k+1}}\times \sum_{j=1}^n M_{nj} (q_{j_{k+1}}-q_{j_{k}})) - \frac{h}{2}q_{n_{k+1}}\times \deriv{V_{k+1}}{q_{n_{k+1}}}
    \end{bmatrix}.\label{eqn:DELw2m}
\end{align}
For a given $(q_{i_k},\omega_{i_k})$, we solve \refeqn{DELw1} to obtain $F_{i_k}\in\SO$. Then, $q_{i_{k+1}}$ and $\omega_{i_{k+1}}$ are computed by \refeqn{qikpw} and \refeqn{DELw2m}, respectively. This yields a discrete flow map in terms of the angular velocity  $(q_{i_k},\omega_{i_k})\mapsto(q_{i_{k+1}},\omega_{i_{k+1}})$, and this process is repeated.
\end{cor}

\subsection{Computational approach}\label{sec:com}
For the discrete equations of motion, we need to solve \refeqn{DEL} and \refeqn{DELw1} to obtain $F_{i_k}\in\SO$. Here we present a computational approach. The implicit equations given by \refeqn{DEL} and \refeqn{DELw1} have the following structure.
\begin{align}
    M_{ii} q_{i}\times F_{i} q_{i}%
    +q_{i}\times \sum_{\substack{j=1\\j\neq i}}^n M_{ij} (F_{j}-I_{3\times 3}) q_{j}=d_i\label{eqn:imp}
\end{align}
for $i\in\{1,\ldots,n\}$, where $M_{ij}\in\Re$, $q_i\in\Sph^2$, $d_i\in\Re^3$ are known, and we need to find $F_i\in\SO$. We derive an equivalent equation in terms of local coordinates for $F_i$. This is reasonable since $F_i$ represents the relative update between two integration steps. Using the Cayley transformation~\cite{Shu.JAS93}, $F_i\in\SO$ can be expressed in terms of $f_i\in\Re^3$ as
\begin{align*}
    F_i & = (I_{3\times 3}+\hat f_i)(I_{3\times 3}-\hat f_i)^{-1},\\
        & = \frac{1}{1+f_i\cdot f_i}((1-f_i\cdot f_i)I_{3\times 3} + 2f_if_i^T +2\hat f_i).
\end{align*}
The operation $F_i q_i$ can be considered as a rotation of the vector $q_i$ about the direction $f_i$ with rotation angle $2\tan^{-1}\norm{f_i}$.  Since the rotation of the vector $q_i$ about the direction $q_i$ has no effect, we can assume that $f_i$ is orthogonal to  $q_i$, i.e. $f_i\cdot q_i=0$. Under this assumption, $F_iq_i$ is given by
\begin{align}
    F_i q_i= \frac{1}{1+f_i\cdot f_i}((1-f_i\cdot f_i)q_i+2\hat f_iq_i).\label{eqn:Fiqi}
\end{align}
Thus, we obtain
\begin{gather*}
    q_i \times F_i q_i = \frac{2}{1+f_i\cdot f_i} q_i\times (f_i\times q_i) = \frac{2}{1+f_i\cdot f_i} f_i,\\
    (F_j-I_{3\times 3})q_j = - \frac{2}{1+f_j\cdot f_j}( q_jf_j^T + \hat q_j )f_j,
\end{gather*}
where we use the property, $\hat q_i f_i = q_i\times f_i=-\hat f_i q_i$. Substituting these into \refeqn{imp}, we obtain
\begin{align}
    \begin{bmatrix} \frac{2M_{11}I_{3\times 3}}{1+f_1\cdot f_1}
       & -\frac{2M_{12}\hat q_1 (\hat q_2+q_2 f_2^T)}{1+f_2\cdot f_2}%
       & \cdots%
       & -\frac{2M_{1n}\hat q_1 (\hat q_n+q_n f_n^T)}{1+f_n\cdot f_n}\\
       -\frac{2M_{21}\hat q_2 (\hat q_1+q_1 f_1^T)}{1+f_1\cdot f_1}%
       & \frac{2M_{22}I_{3\times 3}}{1+f_2\cdot f_2}
       & \cdots
       & -\frac{2M_{2n}\hat q_2 (\hat q_n+q_n f_n^T)}{1+f_n\cdot f_n}\\
       \vdots & \vdots & &\vdots\\
       -\frac{2M_{n1}\hat q_n (\hat q_1+q_1 f_1^T)}{1+f_1\cdot f_1}%
       & -\frac{2M_{n2}\hat q_n (\hat q_2+q_2 f_2^T)}{1+f_2\cdot f_2}%
       & \cdots
       &\frac{2M_{nn}I_{3\times 3}}{1+f_n\cdot f_n}
    \end{bmatrix}
    \begin{bmatrix}    f_1\\f_2\\\vdots\\f_n    \end{bmatrix}
    =    \begin{bmatrix}    d_1\\d_2\\\vdots\\d_n    \end{bmatrix},\label{eqn:impf}
\end{align}
which is an equation equivalent to \refeqn{imp}, written in terms of local coordinates for $F_i$ using the Cayley transformation. Any numerical method to solve nonlinear equations can be applied to find $f_i$. Then, $F_iq_i$ is computed by using \refeqn{Fiqi}. In particular, \refeqn{impf} is written in a form that can be readily applied to a fixed point iteration method~\cite{Kel.BK95}.

If there are no coupling terms in the kinetic energy, we can obtain an explicit solution of \refeqn{imp}. When $M_{ij}=0$ for $i\neq j$, \refeqn{impf} reduces to
\begin{align*}
    \frac{2M_{ii}}{1+f_i\cdot f_i} f_i = d_i.
\end{align*}
Using the identity, $\frac{2\tan \theta}{1+\tan^2\theta}=\sin 2\theta$ for any $\theta\in\Re$, it can be shown that the solution of this equation is given by
$f_i = \tan \parenth{\frac{1}{2}\sin^{-1} (\norm{d_i}/M_{ii})} \frac{d_i}{\norm{d_i}}$.
Substituting this into \refeqn{Fiqi} and rearranging, we obtain
\begin{align*}
    F_i q_i = \frac{d_i}{M_{ii}} \times q_i + \parenth{1-\norm{\frac{d_i}{M_{ii}}}^2}^{1/2} q_i.
\end{align*}
Using this expression, we can rewrite the discrete equations of motion given in \refeqn{DELw1}--\refeqn{DELw2m} in an explicit form.

\begin{cor}
Consider a mechanical system on $(\Sph^2)^n$ whose Lagrangian is expressed as \refeqn{L} where $M_{ij}=0$ for $i\neq j$, i.e. the dynamics are coupled only though the potential energy. The explicit discrete equations of motion are given by
\begin{gather}
    q_{i_{k+1}} = \parenth{h\omega_{i_k} - \frac{h^2}{2M_{ii}} q_{i_k} \times \deriv{V_k}{q_{i_k}}}\times q_{i_k} + \parenth{1-\norm{h\omega_{i_k} - \frac{h^2}{2M_{ii}} q_{i_k} \times \deriv{V_k}{q_{i_k}}}^2}^{1/2} q_{i_k},\label{eqn:qikpexp}\\
    \omega_{i_{k+1}} = \omega_{i_k} - \frac{h}{2M_{ii}} q_{i_k} \times \deriv{V_k}{q_{i_k}}-    \frac{h}{2M_{ii}} q_{i_{k+1}} \times \deriv{V_{k+1}}{q_{i_{k+1}}}\label{eqn:wikpexp}
\end{gather}
for $i\in\{1,\ldots,n\}$.
\end{cor}

\subsection{Properties of variational integrators on $(\Sph^2)^n$}

Since variational integrators are derived from the discrete Hamilton's principle, they are symplectic, and momentum preserving. The discrete action sum can be considered as a zero-form on $(\Sph^2)^n\times(\Sph^2)^n$ which maps the initial condition of a discrete flow satisfying the discrete Euler-Lagrange equation to the action sum for that trajectory. The simplecticity of the discrete flow follows from the fact the the iterated exterior derivative of any differential form is zero. If the discrete Lagrangian exhibits a  symmetry, the corresponding momentum map is preserved since by symmetry, the variation of the discrete Lagrangian in the symmetry direction is zero, which in combination with the discrete Euler--Lagrange equations, implies a discrete version of Noether's theorem. Detailed proofs for the symplectic property and the momentum preserving property can be found in~\cite{MaWe2001}. The total energy oscillates around its initial value with small bounds on a comparatively short timescale, but there is no tendency for the mean of the oscillation in the total energy to drift (increase or decrease) over exponentially long times~\cite{Hai.ANM94}.

The variational integrators presented in this paper preserve the structure of $(\Sph^2)^n$ without need of local parameterizations, explicit constraints or reprojection. Using the characteristics of the homogeneous manifold, the discrete update map is represented by a group action of $\SO$, and a proper subspace is searched to obtain a compact, possibly explicit, form for the numerical integrator. As a result, the following numerical problems are avoided: (i) local parameterizations yield singularities; (ii) numerical trajectories in the vicinity of a singularity experiences numerical ill-conditioning; (iii) unit length constraints lead to additional computational complexity; (iv) reprojection corrupts the numerical accuracy of trajectories~\cite{Hai.BK00,LewNig.JCAM03}.

It can be shown that these variational integrators have second-order accuracy as the discrete action sum is a second-order approximation of the action integral. Higher-order integrators can be easily constructed by applying a symmetric composition method~\cite{Yos.PLA90}.

\subsection{Numerical examples}\label{sec:ne}

The computational properties of variational integrators on $(\Sph^2)^n$ and explicit Runge-Kutta methods are compared for several mechanical systems taken from variety of scientific areas, namely a double spherical pendulum, an $n$-body problem on a sphere, an interconnected system of spherical pendula, pure bending of an elastic rod, a spatial array of magnetic dipoles, and molecular dynamics that evolves on a sphere.

\begin{example}[\textbf{Double Spherical Pendulum}]\label{ex:dsp}
A double spherical pendulum is defined by  two mass particles serially connected to frictionless two degree-of-freedom pivots by rigid massless links acting under a uniform gravitational potential. The dynamics of a double spherical pendulum has been studied in~\cite{MarSchWen.ICIAM95}, and a variational integrator is developed in~\cite{WedMar.PhyD97} by explicitly using unit length constraints.

Let the mass and the length of the pendulum be $m_1,m_2,l_1,l_2\in\Re$, respectively, and let $e_3=[0,0,1]\in\Re^3$ be the direction of gravity. The vector $q_1\in\Sph^2$ represents the direction from the pivot to the first mass, and the vector $q_2\in\Sph^2$ represents the direction from the first mass to the second mass. The inertia matrix is given by $M_{11}=(m_1+m_2)l_1^2$, $M_{12}=m_2l_1l_2$, and $M_{22}=m_2l_2^2$. The gravitational potential is written as $V(q_1,q_2)= -(m_1+m_2)gl_1 e_3\cdot q_1 - m_2 g l_2 e_3 \cdot q_2$ for the gravitational acceleration $g\in\Re$. Substituting these into \refeqn{ELq}--\refeqn{ELwm}, the continuous equations of motion for the double spherical pendulum are given by
\begin{gather}
    \dot q_1 = \omega_1\times q_1\,\quad \dot q_2 = \omega_2\times q_2,\label{eqn:dspqdot}\\
    \begin{bmatrix}%
    (m_1+m_2)l_1^2I_{3\times 3} & -m_2l_1l_2 \hat q_1 \hat q_2\\%
    -m_2l_1l_2 \hat q_2\hat q_1 & m_2l_2^2 I_{3\times 3}
    \end{bmatrix}%
    \begin{bmatrix}
    \dot \omega_1 \\ \dot \omega_2
    \end{bmatrix}
    =
    \begin{bmatrix}
    m_2l_1l_2(\omega_2\cdot\omega_2)\hat q_1 q_2 +(m_1+m_2)gl_1\hat q_1 e_3\\
    m_2l_1l_2(\omega_1\cdot\omega_1)\hat q_2 q_1 + m_2 g l_2 \hat q_2e_3
    \end{bmatrix},\label{eqn:dspwdot}
\end{gather}
which are more compact than existing equations written in terms of angles. Another nice property is that the same structure for the equations of motion is maintained for $n >2$. Thus, it is easy to generalize these equations of motion to a triple, or more generally, a multiple-link spherical pendulum.

\begin{figure}
\centerline{
\subfigure[Trajectory of pendulum]%
    {\includegraphics[width=0.24\textwidth]{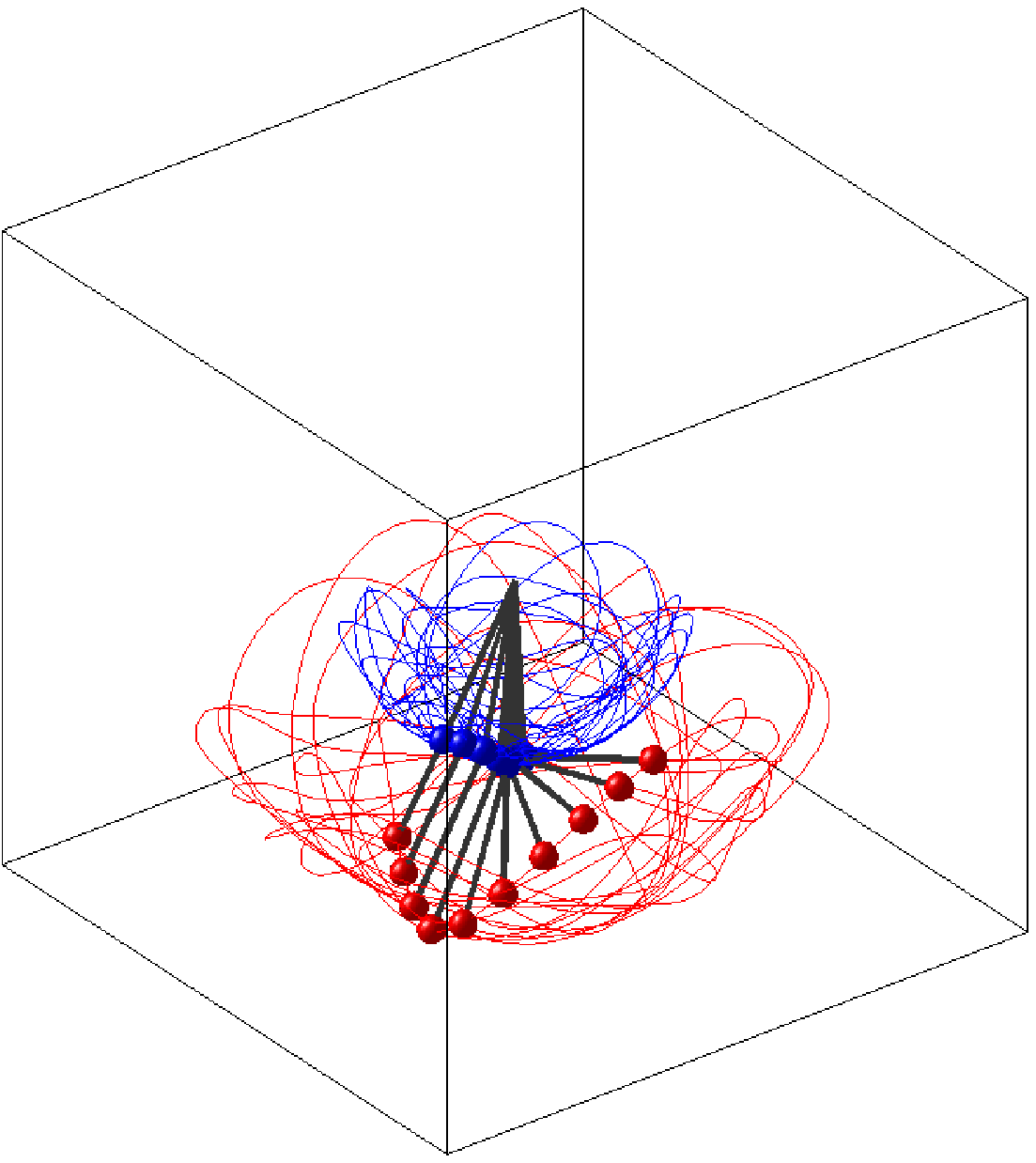}}
\hspace{0.02\textwidth}
\subfigure[Computed total energy]%
    {\includegraphics[width=0.34\textwidth]{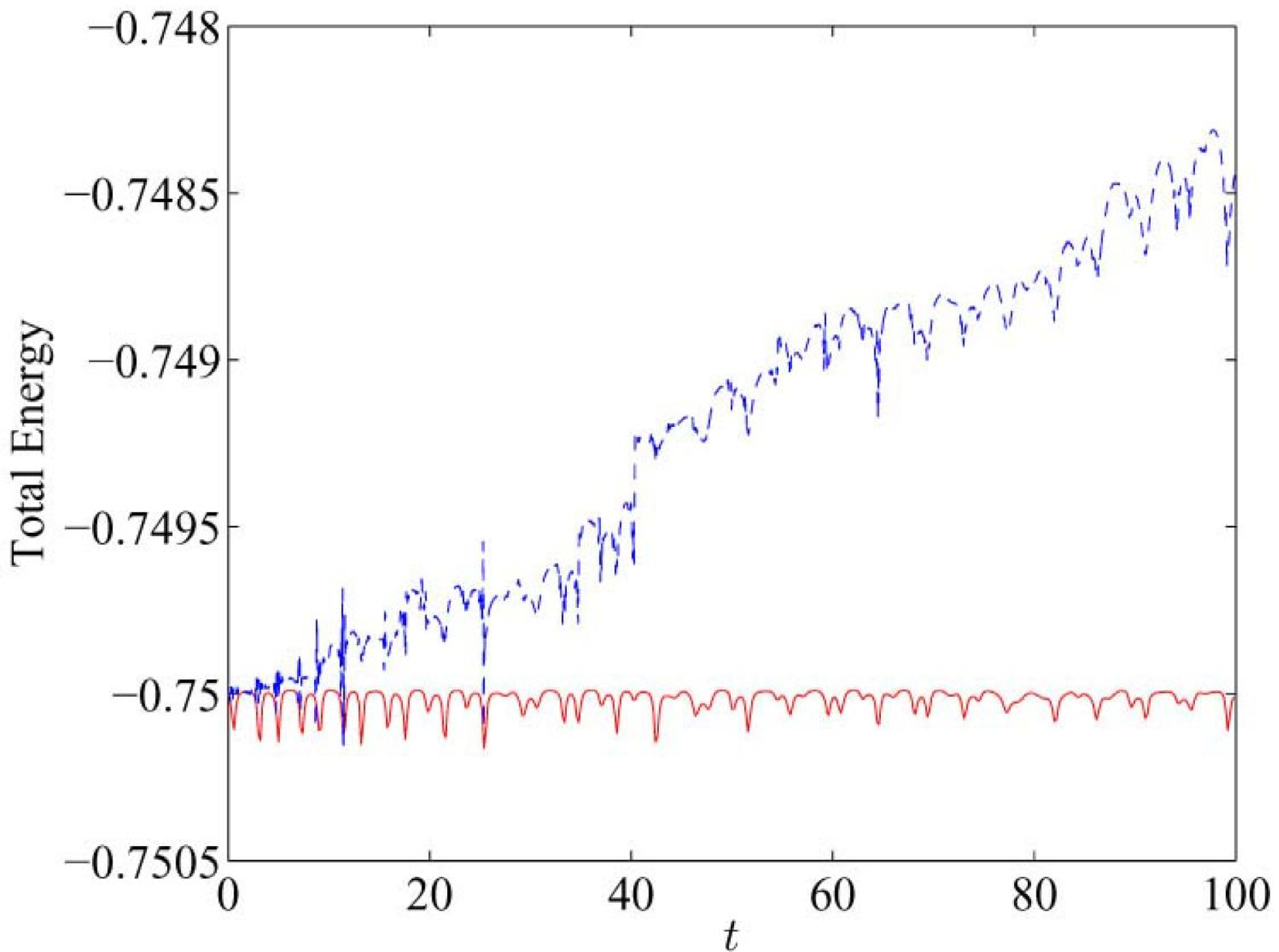}}
\subfigure[Unit length error]%
    {\includegraphics[width=0.34\textwidth]{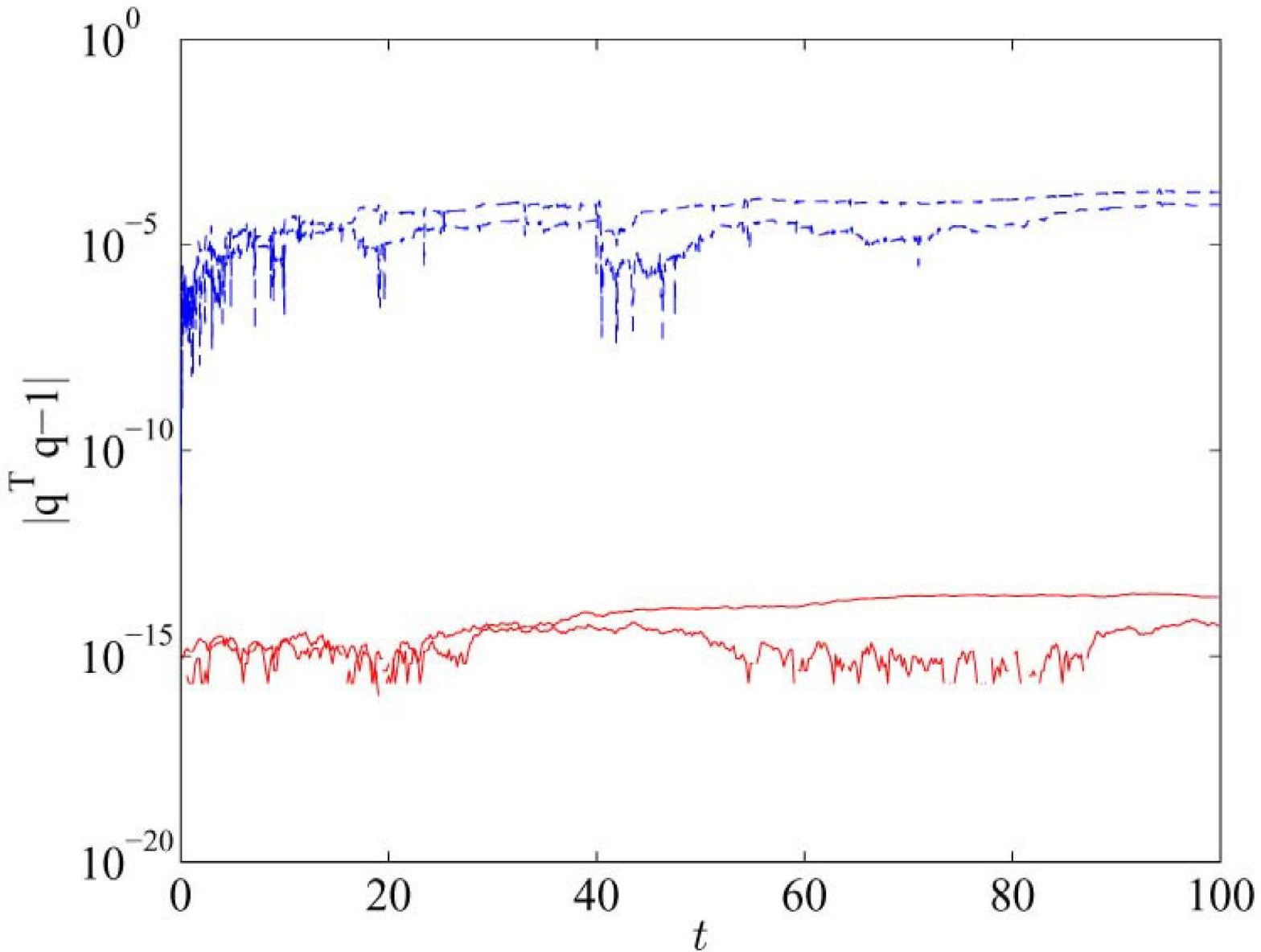}}}
\caption{Numerical simulation of a double spherical pendulum (RK45: blue, dotted, VI: red, solid)}\label{fig:dsp}
\end{figure}

We compare the computational properties of the discrete equations of motion given by \refeqn{DELw1}--\refeqn{DELw2m} with a 4(5)-th order variable step size Runge-Kutta method for \refeqn{dspqdot}--\refeqn{dspwdot}. We choose $m_1=m_2=1\,\mathrm{kg}$, $l_1=l_2=9.81\,\mathrm{m}$. The initial conditions are $q_{1_0}=[0.8660,\,0,\,0.5]$, $q_{2_0}=[0,\,0,\,1]$, $\omega_{1_0}=[-0.4330,\,0,\,0.75]$, $\omega_{2_0}=[0,\,1,\,0]\,\mathrm{rad/sec}$. The simulation time is $100\,\mathrm{sec}$, and the step-size of the discrete equations of motion is $h=0.01$.
\reffig{dsp} shows the computed total energy and the configuration manifold errors. The variational integrator preserves the total energy and the structure of $(\Sph^2)^n$ well for this chaotic motion of the double spherical pendulum. The mean total energy variation is $2.1641\times 10^{-5}\,\mathrm{Nm}$, and the mean unit length error is $8.8893\times 10^{-15}$. But, there is a notable increase of the computed total energy for the Runge-Kutta method, where the mean variation of the total energy is $7.8586\times 10^{-4}\,\mathrm{Nm}$. The Runge-Kutta method also fails to preserve the structure of $(\Sph^2)^n$. The mean unit length error is $6.2742\times 10^{-5}$.
\end{example}

\begin{example}[\textbf{$n$-body Problem on Sphere}]
An $n$-body problem on the two-sphere deals with the motion of $n$ mass particles constrained to lie on a two-sphere, acting under a mutual potential. Let $m_i\in\Re$ and $q_i\in\Sph^2$ be the mass and the position vector of the $i$-th particle, respectively. The $i,j$-th element of the inertia matrix is $M_{ij}=m_i$ when $i=j$, and $M_{ij}=0$ otherwise. In~\cite{KozHar.CMDA92}, the following expression for the potential is introduced as an analog of a gravitational potential, \begin{align*}
    V (q_1,\ldots, q_n) = -\frac{\gamma}{2}\sum_{\substack{i,j=1\\i\neq j}}^n \frac{q_i\cdot q_j}{\sqrt{1-(q_i\cdot q_j)^2}}
\end{align*}
for a constant $\gamma$. Substituting these into \refeqn{EL}, the continuous equations of motion for the $n$-body problem on a sphere are given by
\begin{align}
    m_i \ddot q_i = -m_i (\dot q_i\cdot \dot q_i)q_i - q_i\times\big(q_i\times\gamma\sum_{\substack{j=1\\j\neq i}}^n \frac{q_j}{(1-(q_i\cdot q_j)^2)^{3/2}}\big)\label{eqn:ELnbpsph}
\end{align}
for $i\in\{1,\ldots,n\}$.

A two-body problem on the two-sphere under this gravitational potential is studied in~\cite{HaLuWa2003} by explicitly using unit length constraints. Here we study a three-body problem, $n=3$. Since there are no coupling terms in the kinetic energy, we use the explicit form of the variational integrator. We compare the computational properties of the discrete equations of motion given by \refeqn{qikpexp}--\refeqn{wikpexp} with a 2-nd order fixed step size Runge-Kutta method for \refeqn{ELnbpsph}. We choose $m_1=m_2=m_3=1$, and $\gamma=1$. The initial conditions are $q_{1_0}=[0,\,-1,\,0]$, $q_{2_0}=[0,\,0,\,1]$, $q_{3_0}=[-1,\,0,\,0]$, $\omega_{1_0}=[0,\,0,\,-1.1]$, $\omega_{2_0}=[1,\,0,\,0]$, and $\omega_{3_0}=[0,\,1,\,0]$. The simulation time is $10\,\mathrm{sec}$. \reffig{nbpsph} shows the computed total energy and the unit length errors for various step sizes. The total energy variations and the unit length errors for the variational integrator are smaller than those of the Runge-Kutta method for the same time step size by several orders of magnitude. For the variational integrator, the total energy error is reduced by almost 100 times from $1.1717\times 10^{-4}$ to $1.1986\times 10^{-6}$ when the step size is reduced by 10 times from $10^{-3}$ to $10^{-4}$, which verifies the second order accuracy numerically.

\begin{figure}
\centerline{
    \subfigure[Trajectory of particles]{%
        \includegraphics[width=0.26\textwidth]{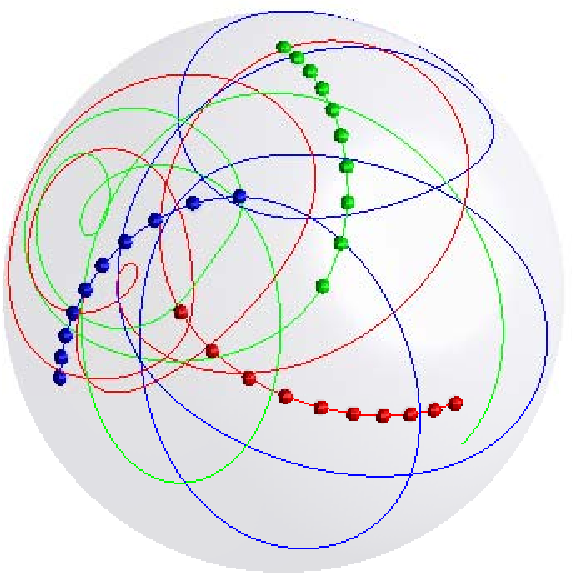}}
    \hspace{0.02\textwidth}
    \subfigure[Total energy error v.s. step size]{%
        \includegraphics[width=0.34\textwidth]{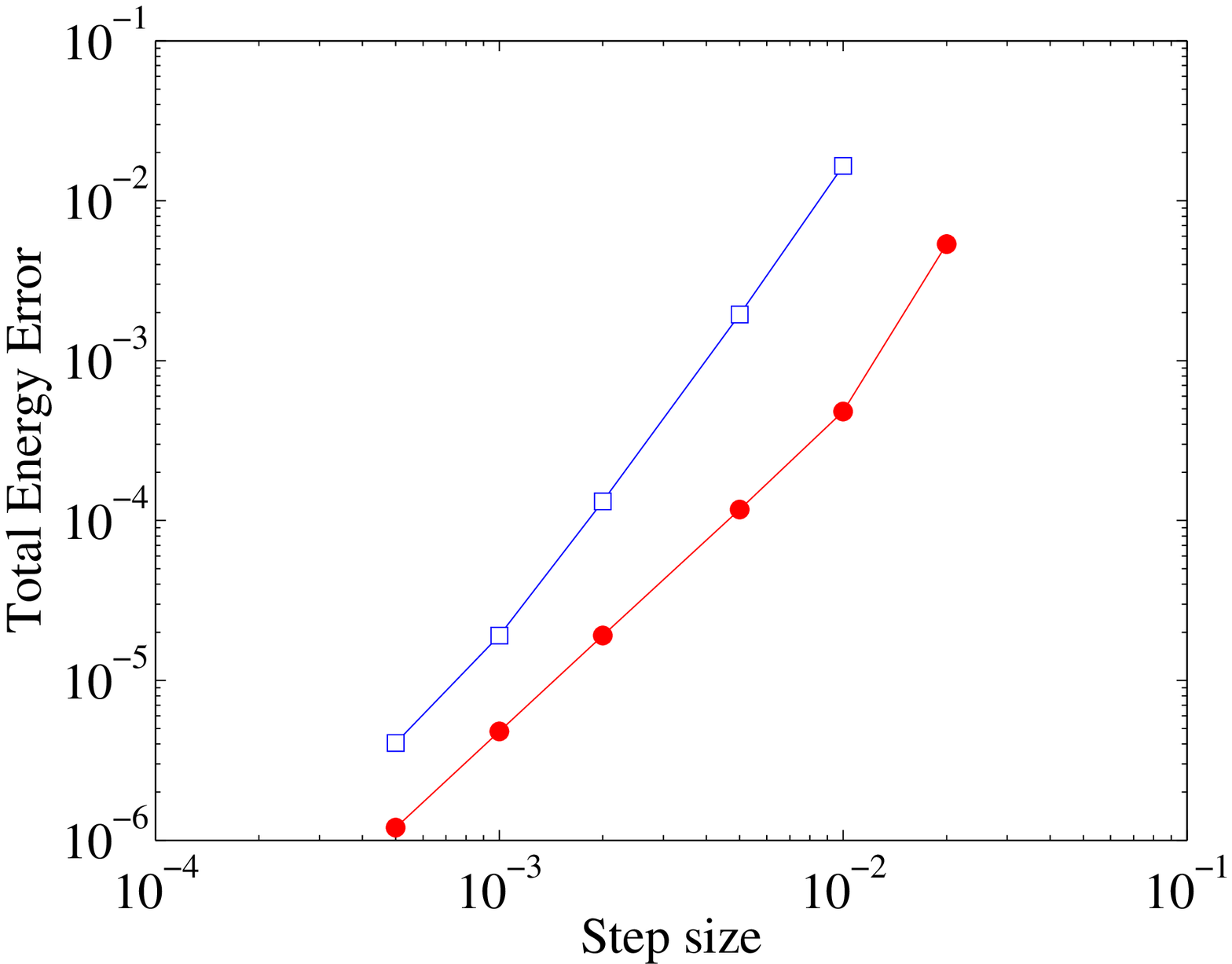}}
    \subfigure[Unit length error v.s. step size]{%
        \includegraphics[width=0.34\textwidth]{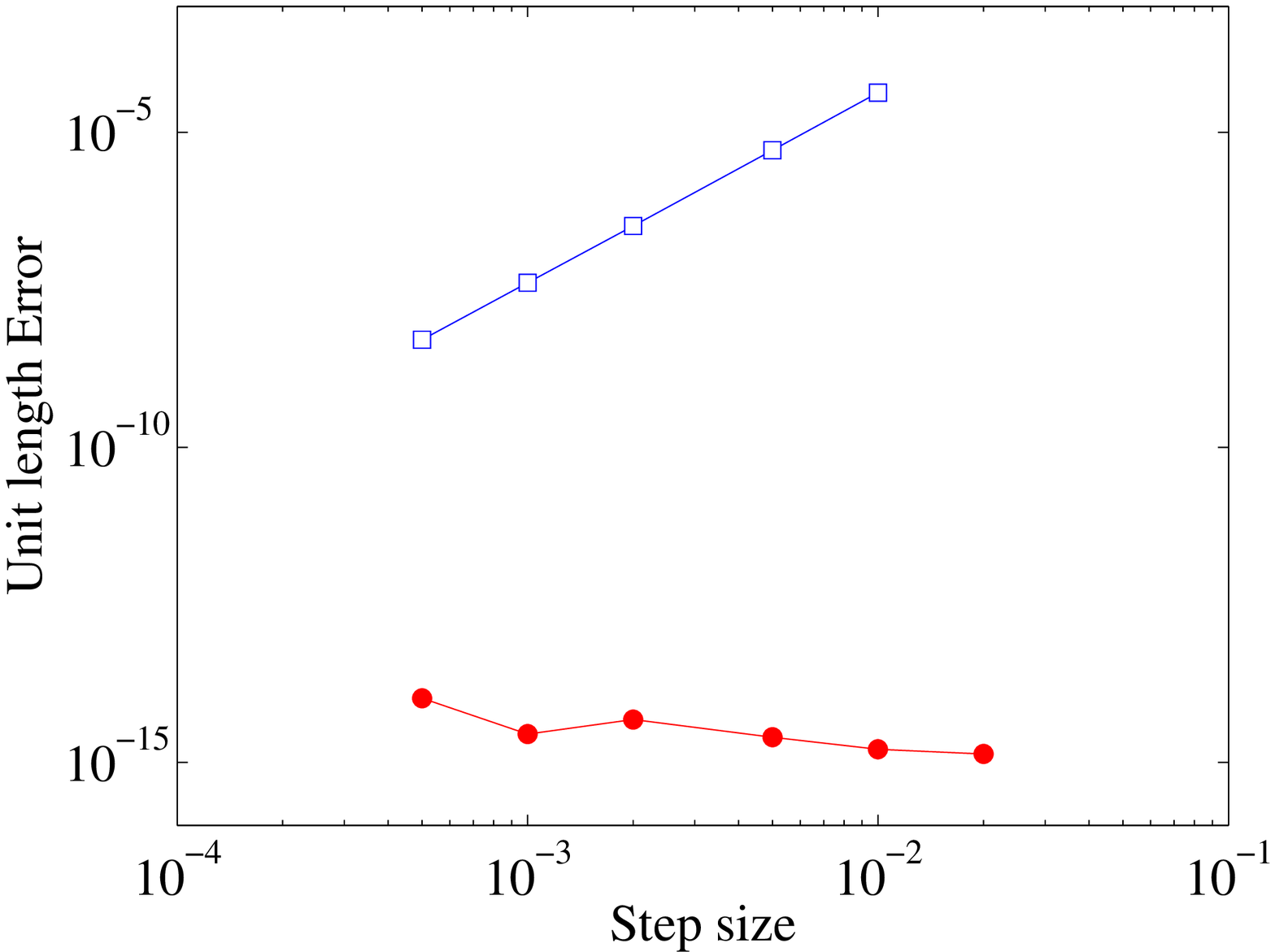}}}
\caption{Numerical simulation of a 3-body problem on sphere (RK2: blue, square, VI: red, circle)}\label{fig:nbpsph}
\end{figure}
\end{example}

\begin{example}[\textbf{Interconnection of Spherical Pendula}]
We study the dynamics of $n$ spherical pendula connected by linear springs. Each pendulum is a mass particle connected to a frictionless two degree-of-freedom pivot by a rigid massless link acting under a uniform gravitational potential. It is assumed that all of the pivot points lie on a horizontal plane, and some pairs of pendulua are connected by linear springs at the centers of links.

Let the mass and the length of the $i$-th pendulum be $m_i,l_i\in\Re$, respectively. The vector $q_i\in\Sph^2$ represents the direction from the $i$-th pivot to the $i$-th mass. The inertia matrix is given by $M_{ij}=m_il_i^2$ when $i=j$, and $M_{ij}=0$ otherwise. Let $\Xi$ be a set defined such that $(i,j)\in\Xi$ if the $i$-th pendulum and the $j$-th pendulum are connected. For a connected pair $(i,j)\in\Xi$, define $\kappa_{ij}\in\Re$ and $r_{ij}\in\Re^3$ as the corresponding spring constant and the vector from the $i$-th pivot to the $j$-th pivot, respectively. The bases for the inertial frame are chosen such that the direction along  gravity is denoted by $e_3=[0,0,1]\in\Re^3$, and the horizontal plane is spanned by  $e_1=[0,0,1],e_2=[0,1,0]\in\Re^3$. The potential energy is given by
\begin{align*}
    V(q_1,\ldots q_n) = -\sum_{i=1}^n m_igl_iq_i\cdot e_3%
     + \sum_{(i,j)\in\Xi}\frac{1}{2}\kappa_{ij}
     \parenth{\norm{r_{ij}+\frac{1}{2}l_jq_j-\frac{1}{2}l_iq_i}-\norm{r_{ij}}}^2.
\end{align*}
Substituting these into \refeqn{ELw}--\refeqn{ELq}, the continuous equations of motion for the interconnection of spherical pendula are given by
\begin{gather}
    m_il_i^2 \dot\omega_i = - q_i\times\deriv{V}{q_i}\label{eqn:ELwsyspend}\\
    \dot q_i = \omega_i\times q_i\label{eqn:ELqsyspend}
\end{gather}
for $i\in\{1,\ldots,n\}$.

We compare the computational properties of the discrete equations of motion given by \refeqn{qikpexp}--\refeqn{wikpexp} with a 2-nd order fixed step size explicit Runge-Kutta method for \refeqn{ELwsyspend}--\refeqn{ELqsyspend}, and the same Runge-Kutta method with reprojection; at each time step, the vectors $q_{i_k}$ are projected onto $\Sph^2$ by using normalization.

We choose four interconnected pendula, $n=4$, and we assume each pendulum has the same mass and length; $m_i=0.1\,\mathrm{kg}$, $l_i=0.1\,\mathrm{m}$. The pendula are connected as $\Xi=\{(1,2),(2,3),(3,4),(4,1)\}$, and the corresponding spring constants and the relative vector between pivots are given by $\kappa_{12}=10$, $\kappa_{12}=20$, $\kappa_{12}=30$, $\kappa_{12}=40\,\mathrm{N/m}$, $r_{12}=-r_{34}=l_ie_1$, and $r_{23}=-r_{41}=-l_ie_2$. The initial conditions are chosen as $q_{1_0}=q_{2_0}=q_{4_0}=e_3$, $q_{3_0}=[0.4698,0.1710,0.8660]$, $\omega_{1_0}=[-10,4,0]$, and $\omega_{2_0}=\omega_{3_0}=\omega_{4_0}=0\,\mathrm{rad/sec}$

\begin{figure}
\centerline{
    \subfigure[Motion of pendula]{%
        \includegraphics[width=0.27\textwidth]{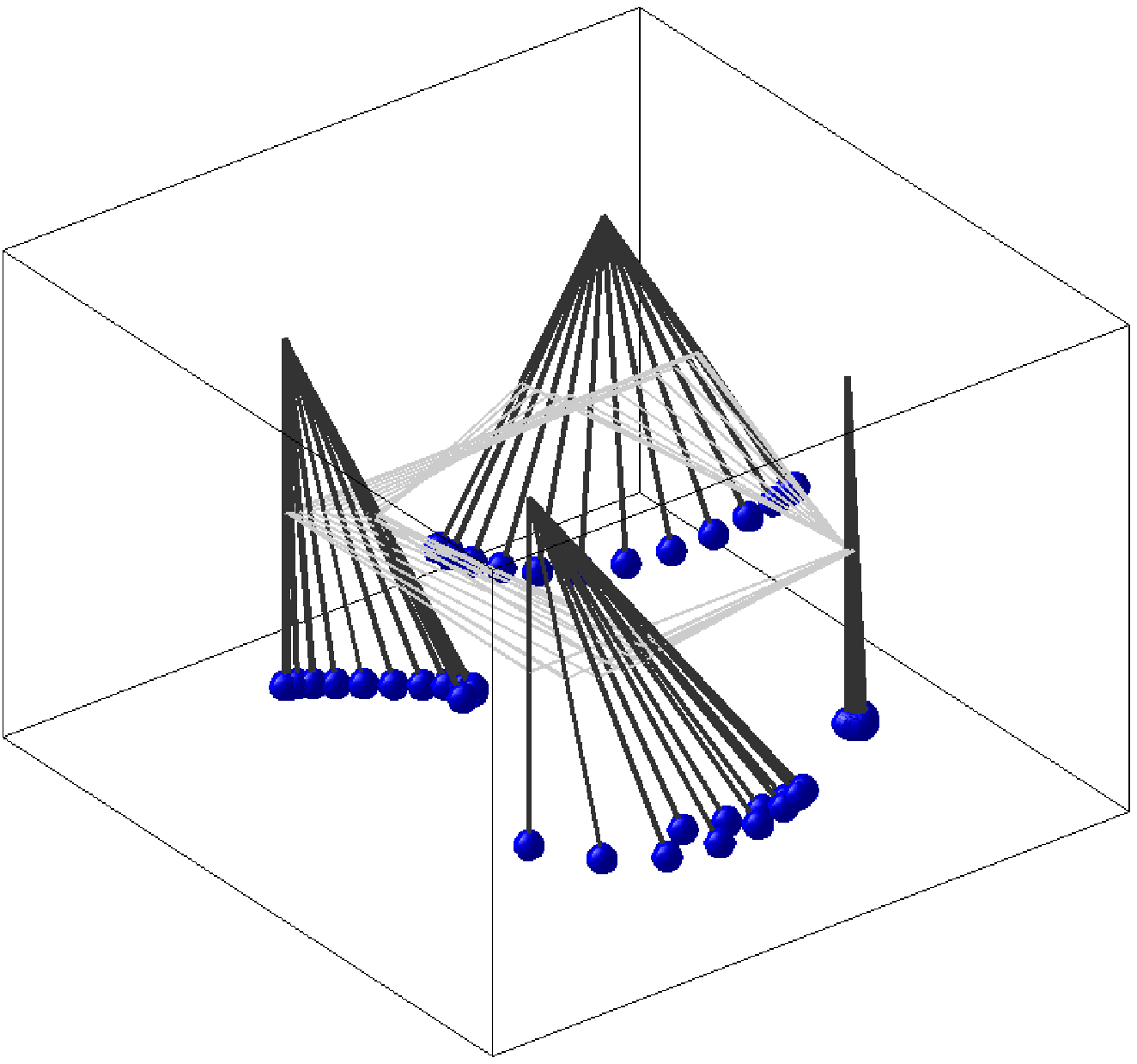}}
    \hspace{0.02\textwidth}
    \subfigure[Computed total energy]{%
        \includegraphics[width=0.34\textwidth]{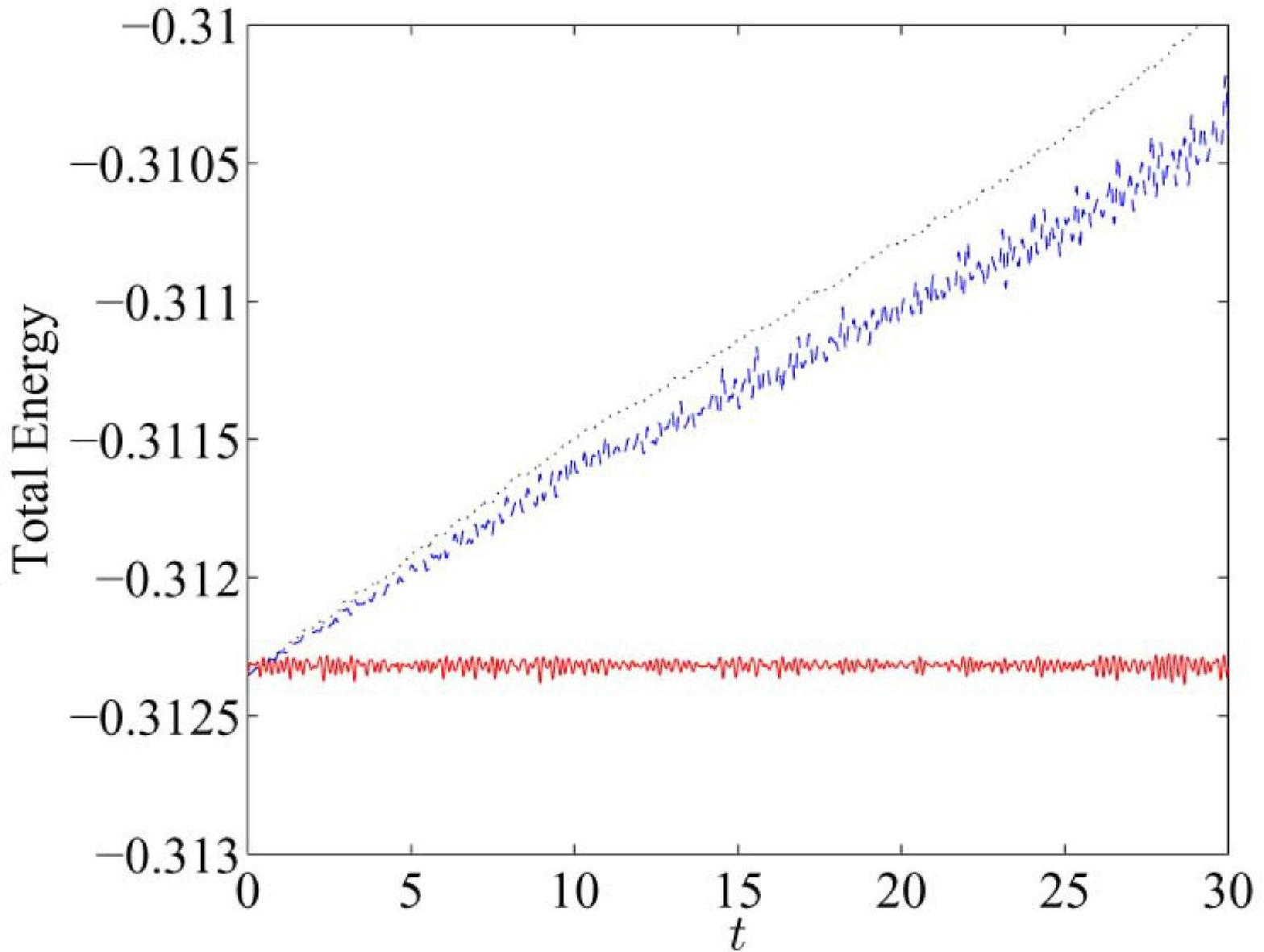}}
    \subfigure[Unit length error]{%
        \includegraphics[width=0.34\textwidth]{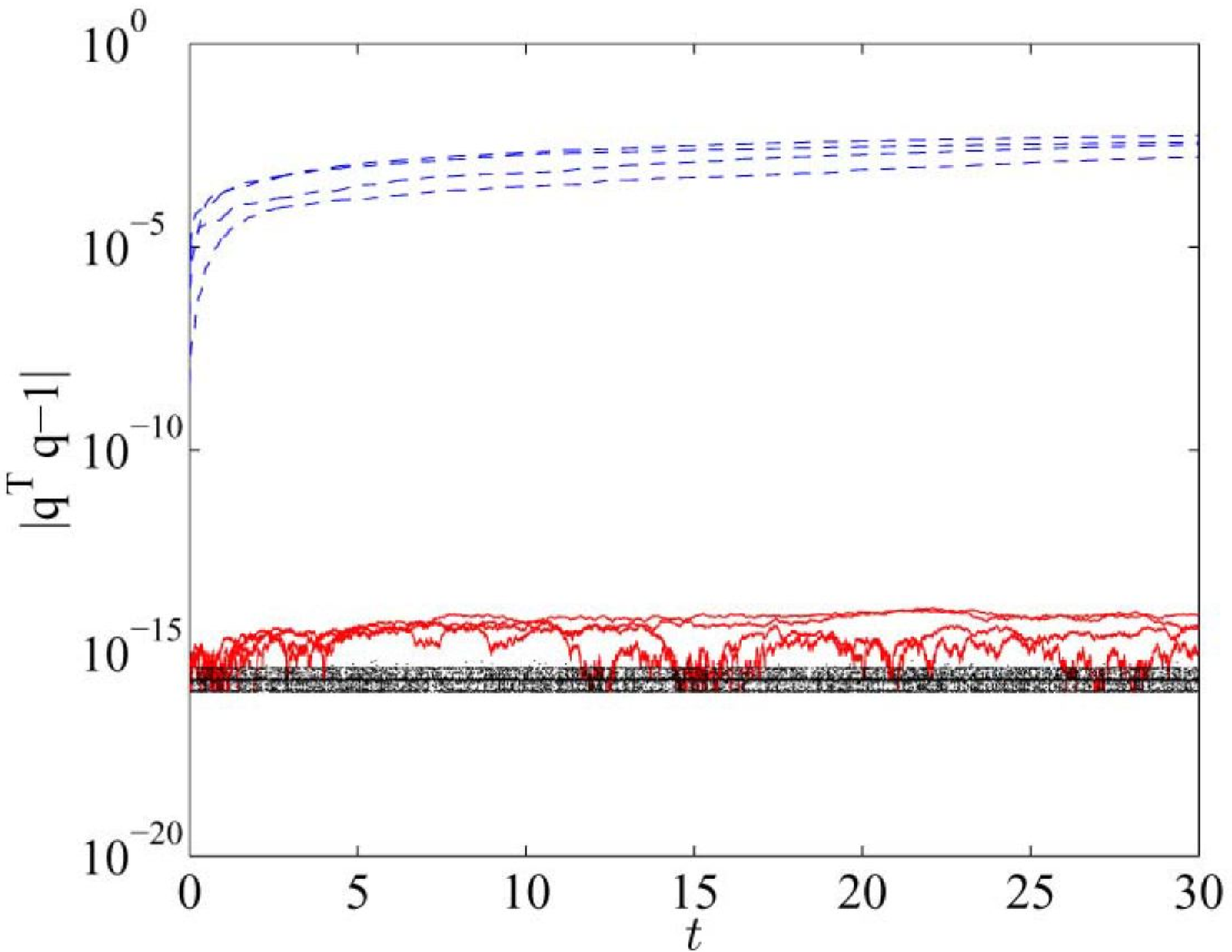}}}
\caption{Numerical simulation of a system of 4 spherical pendula (RK2: blue, dotted, RK2 with projection: black, dashed, VI: red, solid)}\label{fig:syspend}
\end{figure}

\reffig{syspend} shows the computed total energy and the unit length errors. The variational integrator preserves the total energy and the structure of $(\Sph^2)^n$ well. The mean total energy variation is $3.6171\times 10^{-5}\,\mathrm{Nm}$, and the mean unit length error is $4.2712\times 10^{-15}$. For both Runge-Kutta methods, there is a notable increase of the computed total energy. It is interesting to see that the reprojection approach makes the total energy error worse, even though it preserves the structure of $(\Sph^2)^n$ accurately. This shows that a standard reprojection method can corrupt numerical trajectories~\cite{Hai.BK00,LewNig.JCAM03}.
\end{example}

\begin{example}[\textbf{Pure Bending of Elastic Rod}]\label{ex:rod}
We study the dynamics of $(n+1)$ rigid rod elements that are serially connected by rotational springs, where the `zeroth' rod is assumed to be fixed to a wall. Thus, the configuration space is $(\Sph^2)^n$. This can be considered as a simplified dynamics model for pure non-planar bending of a thin elastic rod that is clamped at one end and free at the other end. Notably, this approach is geometrically exact, and preserves the length of the elastic rod in the presence of large displacements.

The mass and the length of the $i$-th rod element are denoted by $m_i,l_i\in\Re$, respectively. The inertia matrix is given by
\begin{align*}
    M_{ii}=\frac{1}{3}m_il_i^2 + \sum_{k=i+1}^n m_i l_i^2,\quad M_{ij}=\sum_{k=\max\{i,j\}}^n \frac{1}{2} m_k l_k^2
\end{align*}
for $i,j\in\{1,\ldots\,n\}$ and $i\neq j$. The potential energy is composed of gravitational terms and elastic bending terms given by
\begin{align*}
    V(q_1,\ldots,q_n) = -\sum_{i=1}^n m_i g \big(\sum_{j=1}^{i-1}l_jq_j+\frac{1}{2}l_iq_i\big) \cdot e_3  +\frac{1}{2}\kappa_i (1-q_{i-1}\cdot q_i)^2,
\end{align*}
where a constant vector $q_0\in\Sph^2$ denotes the direction of the zeroth rod element fixed to a wall, and $\kappa_i\in\Re$ denotes spring constants. The bases for the inertial frame are chosen such  that the gravity direction is denoted by $e_3=[0,0,1]\in\Re^3$, and the horizontal plane is spanned by  $e_1=[0,0,1],e_2=[0,1,0]\in\Re^3$.
Suppose that the total mass and length of rod are given by $m, l$, and each rod element has the same mass and length, i.e. $m_i=\frac{m}{n+1}$, $l_i=\frac{l}{n+1}$ for $i\in\{0,\ldots,n\}$. Substituting these into \refeqn{ELm}, the continuous equations of motion for the pure bending of an elastic rod are given by
\begin{align}
    \begin{bmatrix}%
    \frac{n-2/3}{(n+1)^3}ml^2 I_{3\times 3} & -\frac{n-1}{2(n+1)^3}ml^2 \hat q_1 \hat q_1 & \cdots & -\frac{1}{2(n+1)^3}ml^2\hat q_1 \hat q_1\\%
    -\frac{n-1}{2(n+1)^3}ml^2 \hat q_2\hat q_2 & \frac{n-5/3}{(n+1)^3}ml^2 I_{3\times 3} & \cdots & -\frac{1}{2(n+1)^3}ml^2 \hat q_2 \hat q_2\\%
    \vdots & \vdots & & \vdots\\
    -\frac{1}{2(n+1)^3}ml^2 \hat q_n \hat q_n & -\frac{1}{2(n+1)^3}ml^2\hat q_n \hat q_n & \cdots & \frac{1/3}{(n+1)^3}ml^2 I_{3\times 3}
    \end{bmatrix}%
    \begin{bmatrix}
    \ddot q_1 \\ \ddot q_2 \\ \vdots \\ \ddot q_n
    \end{bmatrix}
    =
    \begin{bmatrix}
    -\frac{n-2/3}{(n+1)^3}ml^2(\dot q_1 \cdot \dot q_1) q_1 +\hat q_1^2 \deriv{V}{q_1}\\
    -\frac{n-5/3}{(n+1)^3}ml^2(\dot q_2 \cdot \dot q_2) q_2 +\hat q_2^2 \deriv{V}{q_2}\\
    \vdots\\
    -\frac{1/3}{(n+1)^3}ml^2(\dot q_n \cdot \dot q_n) q_n +\hat q_n^2 \deriv{V}{q_n}
    \end{bmatrix},\label{eqn:ELrod}
\end{align}

We compare the computational properties of the discrete equations of motion given by \refeqn{DELw1}--\refeqn{DELw2m} with a 4(5)-th order variable step size Runge-Kutta method for \refeqn{ELrod}. We choose 10 rod elements, $n=10$, and the total mass and the total length are $m=55\,\mathrm{g}$, $l=1.1,\mathrm{m}$. The spring constants are chosen as $\kappa_i=1000\,\mathrm{Nm}$. Initially, the rod is aligned horizontally; $q_{i_0}=e_1$ for all $i\in{1,\ldots n}$. The initial angular velocity for each rod element is zero except $\omega_{5_0}=[0,0,10]\,\mathrm{rad/sec}$. This represents the  dynamics of the rod after an initial impact. The simulation time is $3\,\mathrm{sec}$, and the step-size of the discrete equations of motion is $h=0.0001$.

\begin{figure}
\centerline{
    \subfigure[Deformation of rod]{%
        \includegraphics[width=0.28\textwidth]{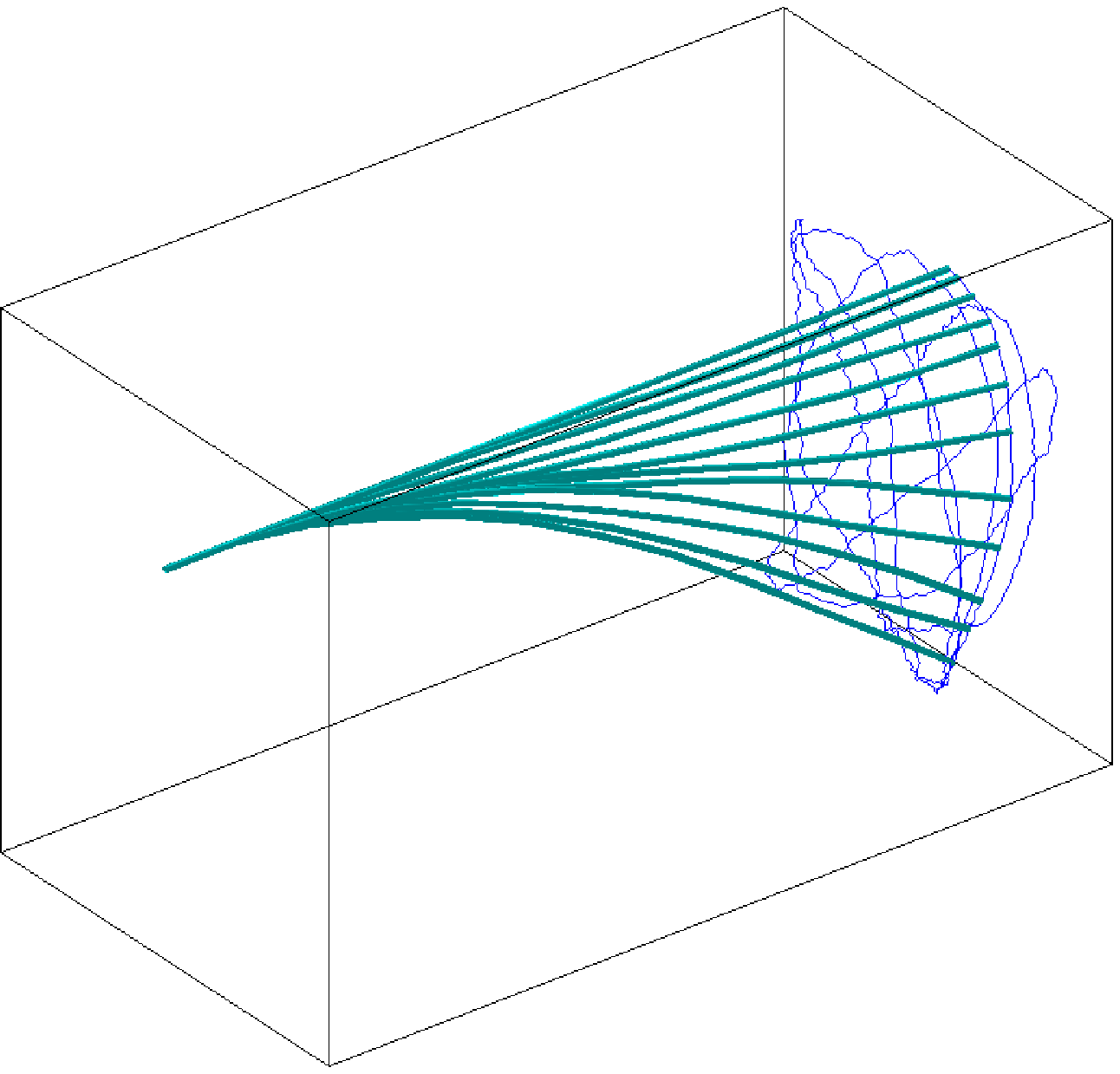}}
    \hspace{0.02\textwidth}
    \subfigure[Computed total energy]{%
        \includegraphics[width=0.34\textwidth]{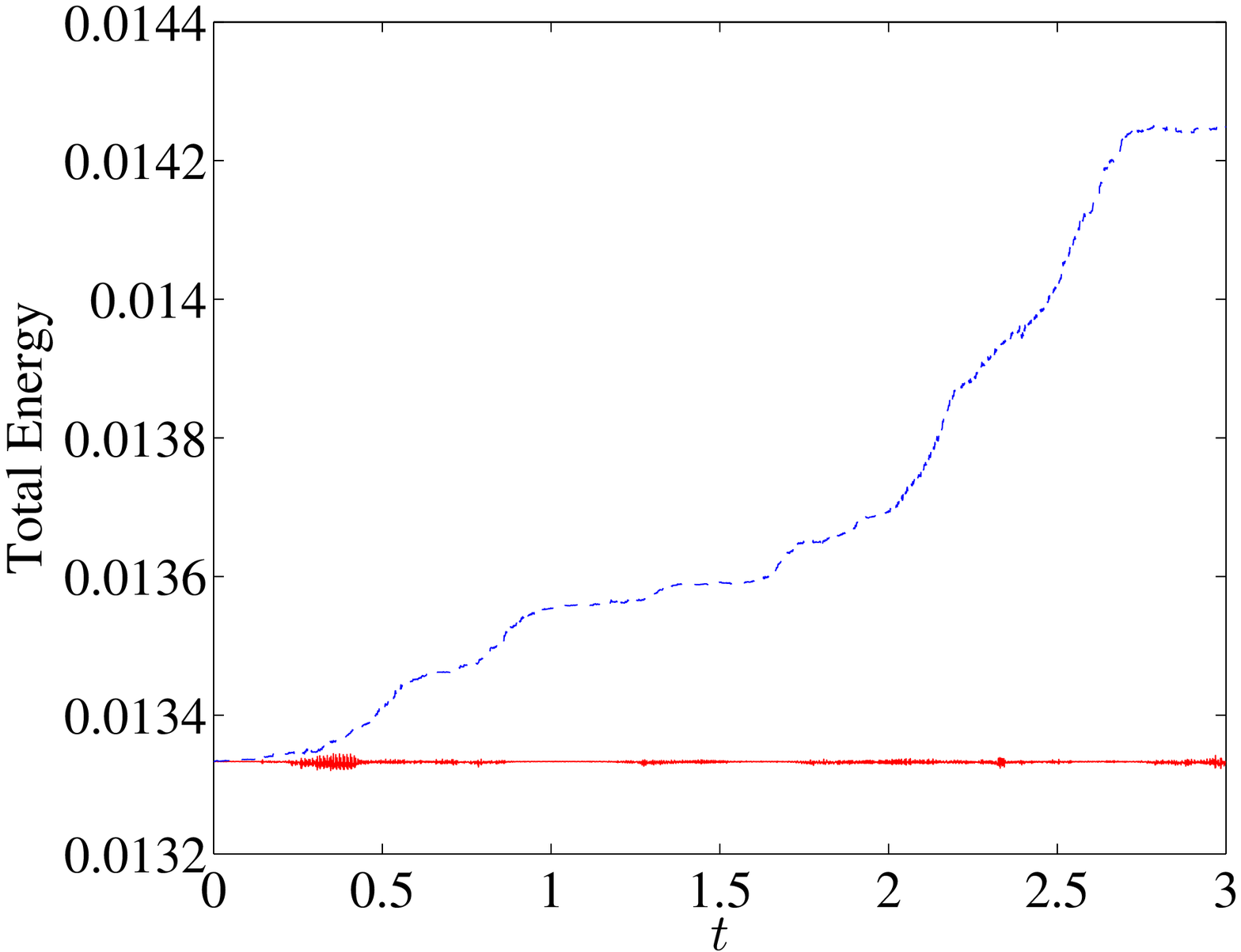}}
    \subfigure[Unit length error]{%
        \includegraphics[width=0.34\textwidth]{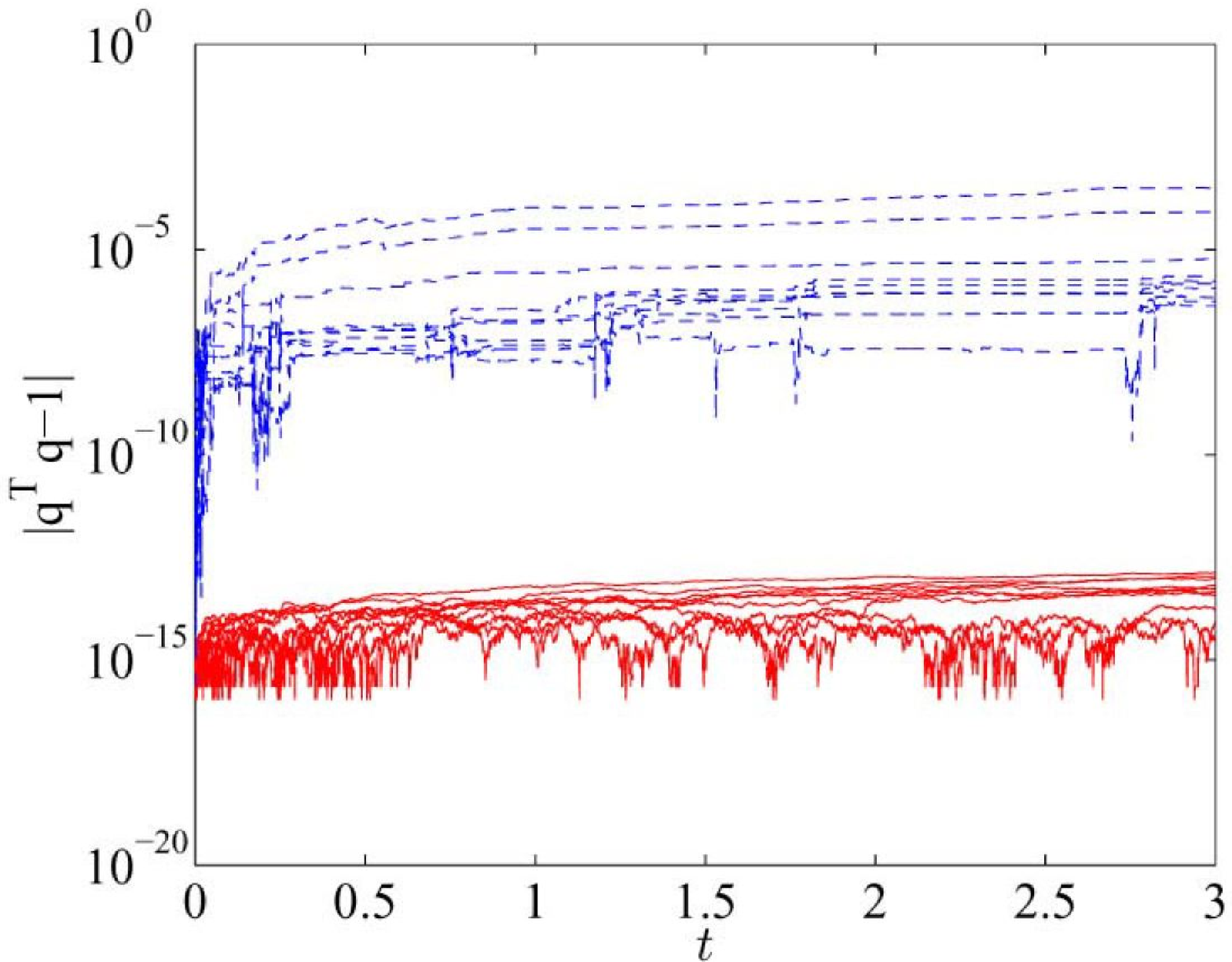}}}
\caption{Numerical simulation of an elastic rod (RK45: blue, dotted, VI: red, solid)}\label{fig:rod}
\end{figure}

\reffig{rod} shows the computed total energy and the unit length errors. The variational integrator preserves the total energy and the structure of $(\Sph^2)^n$. The mean total energy variation is $1.4310\times 10^{-6}\,\mathrm{Nm}$, and the mean unit length error is $2.9747\times 10^{-14}$. There is a notable dissipation of the computed total energy for the Runge-Kutta method, where the mean variation of the total energy is $3.5244\times 10^{-4}\,\mathrm{Nm}$. The Runge-Kutta method also fail to preserve the structure of $(\Sph^2)^n$. The mean unit length error is $1.8725\times 10^{-5}$.
\end{example}

\begin{example}[\textbf{Spatial Array of Magnetic Dipoles}]
We study dynamics of $n$ magnetic dipoles uniformly distributed on a plane. Each magnetic dipole is modeled as a spherical compass; a thin rod magnet supported by a frictionless, two degree-of-freedom pivot acting under their mutual magnetic field. This can be considered as a simplified model for the dynamics of micromagnetic particles~\cite{CheJalLee.PRB06}.

The mass and the length of the $i$-th magnet are denoted by $m_i,l_i\in\Re$, respectively. The magnetic dipole moment of the $i$-th magnet is denoted by $\nu_i q_i$, where $\nu_i\in\Re$ is the constant magnitude of the magnetic moment measured in ampere square-meters, and $q_i\in\Sph^2$ is the direction of the north pole from the pivot point. Thus, the configuration space is $(\Sph^2)^n$. The inertia matrix is given by $M_{ij}=\frac{1}{12}m_il_i^2$ when $i=j$, and $M_{ij}=0$ otherwise. Let $r_{ij}\in\Re^3$ be the vector from the $i$-pivot point to the $j$-th pivot point. The mutual potential energy of the array of magnetic dipoles are given by
\begin{align*}
    V(q_1,\ldots,q_n) = \frac{1}{2}\sum_{\substack{i,j=1\\j\neq i}}^n \frac{\mu\,\nu_i\nu_j}{4\pi \|r_{ij}\|^3}\bracket{(q_i\cdot q_j) - \frac{3}{\|r_{ij}\|^2} ( q_i\cdot r_{ij})( q_j \cdot r_{ij})},
\end{align*}
where $\mu=4\pi\times 10^{-7}\,\mathrm{N\cdot A^{-2}}$ is the permeability constant.
Substituting these into \refeqn{ELw}--\refeqn{ELq}, the continuous equations of motion for the spatial array of magnetic dipoles are given by
\begin{gather}
    \frac{1}{12}m_il_i^2 \dot\omega_i = - q_i\times\sum_{\substack{j=1\\j\neq i}}^n
    \frac{\mu\,\nu_i\nu_j}{4\pi \|r_{ij}\|^3}\bracket{q_j - \frac{3}{\|r_{ij}\|^2} r_{ij}( q_j \cdot r_{ij})},
    \label{eqn:ELwmag}\\
        \dot q_i = \omega_i\times q_i\label{eqn:ELqmag}
\end{gather}
for $i\in\{1,\ldots,n\}$.

We compare the computational properties of the discrete equations of motion given by \refeqn{qikpexp}--\refeqn{wikpexp} with a 4(5)-th order variable step size Runge-Kutta method for \refeqn{ELwmag}--\refeqn{ELqmag}. We choose 16 magnetic dipoles, $n=16$, and we assume each magnetic dipole has the same mass, length, and magnitude of magnetic moment; $m_i=0.05\,\mathrm{kg}$, $l_i=0.02\,\mathrm{m}$, $\nu_i=0.1\,\mathrm{A\cdot m^2}$. The magnetic dipoles are located at vertices of a $4\times 4$ square grid in which the edge of a unit square has the length of $1.2l_i$. The initial conditions are chosen as $q_{i_0}=[1,0,0]$, $\omega_{i_0}=[0,0,0]$ for all $i\in\{1,\ldots,16\}$ except $q_{16_0}=[0.3536, 0.3536, -0.8660]$ and $\omega_{1_0}=[0,0.5,0]\,\mathrm{rad/sec}$.

\begin{figure}
\centerline{
    \subfigure[Motion of magnetic dipoles]{%
        \includegraphics[width=0.36\textwidth]{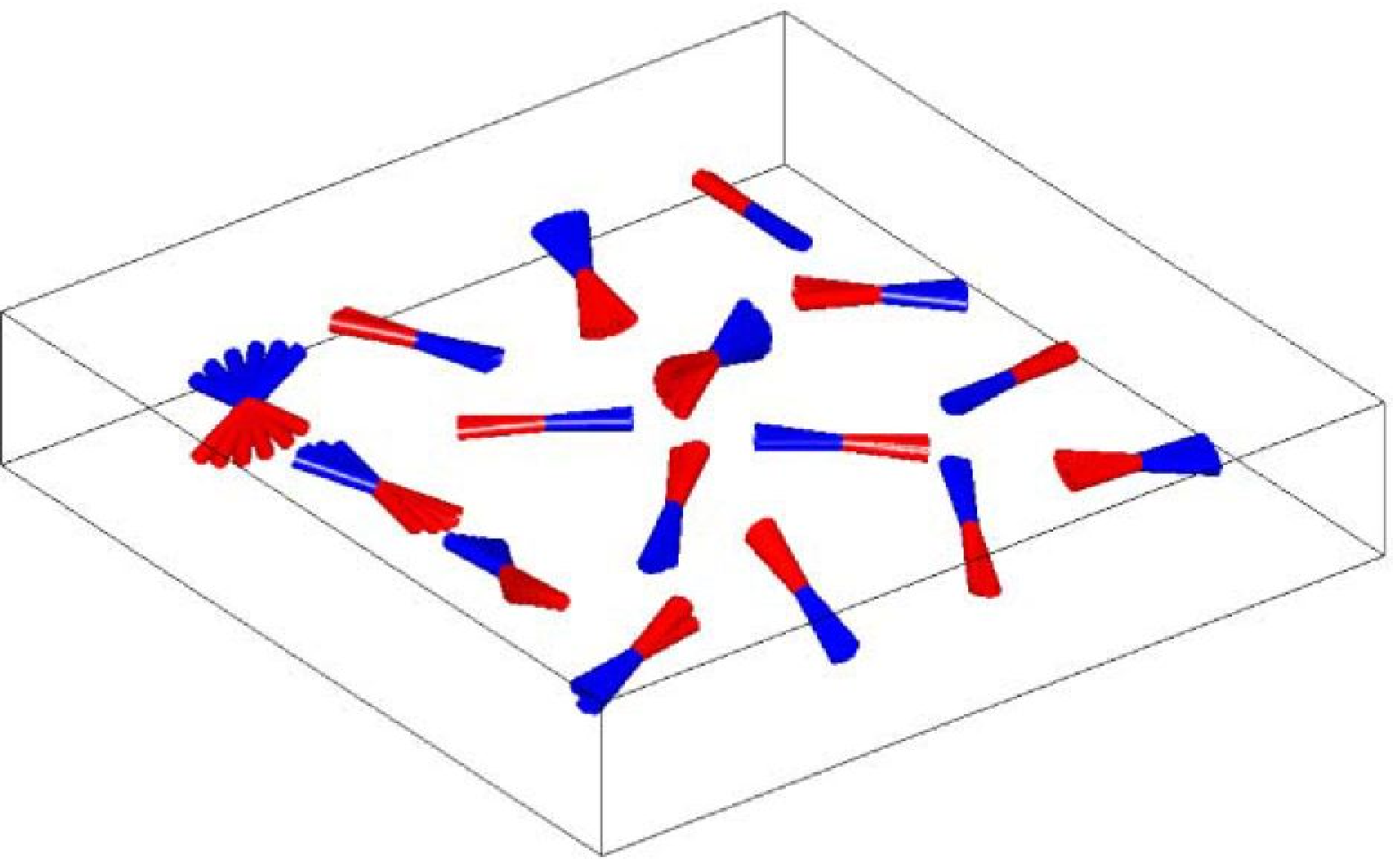}}
    \subfigure[Computed total energy]{%
        \includegraphics[width=0.32\textwidth]{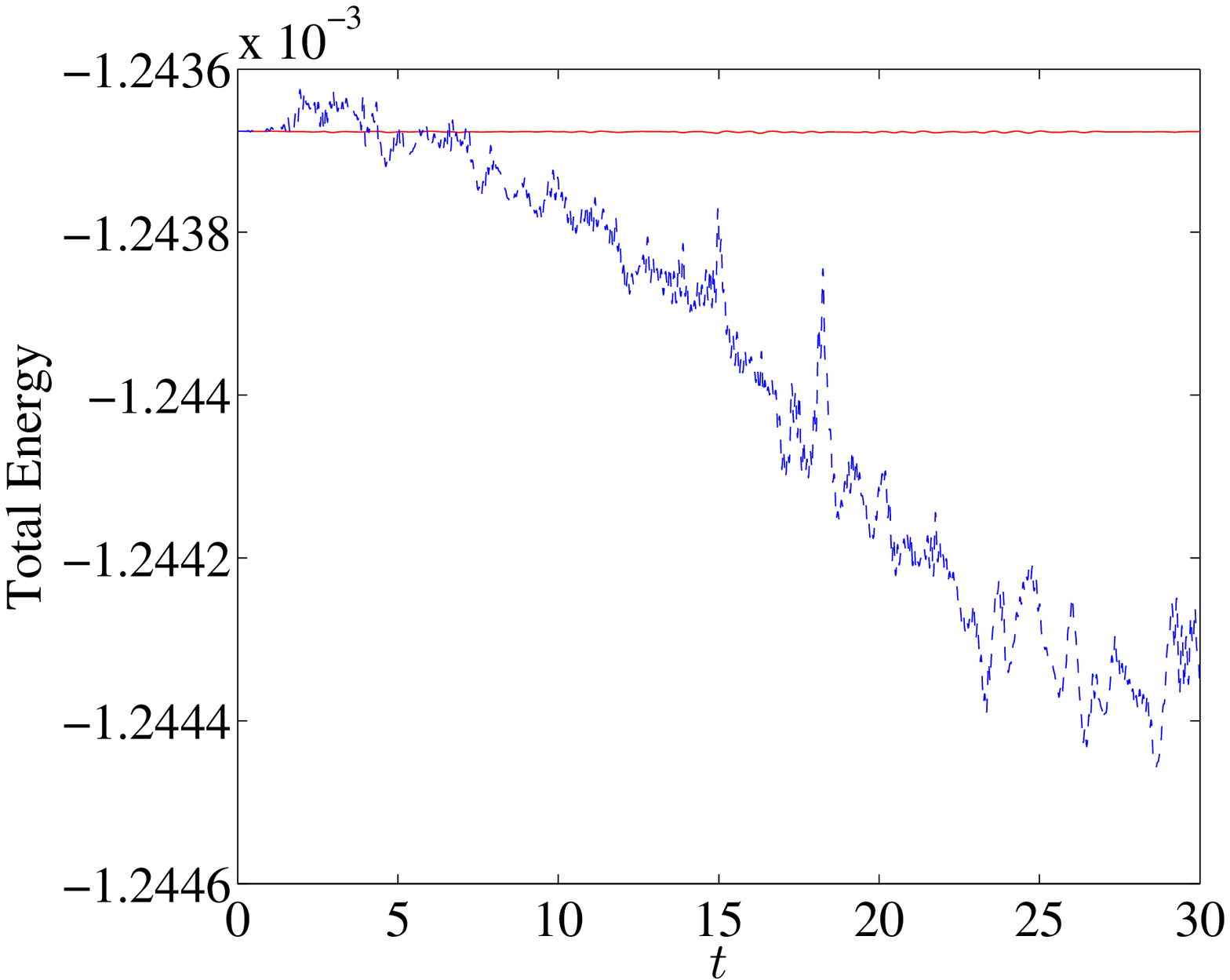}}
    \subfigure[Unit length error]{%
        \includegraphics[width=0.32\textwidth]{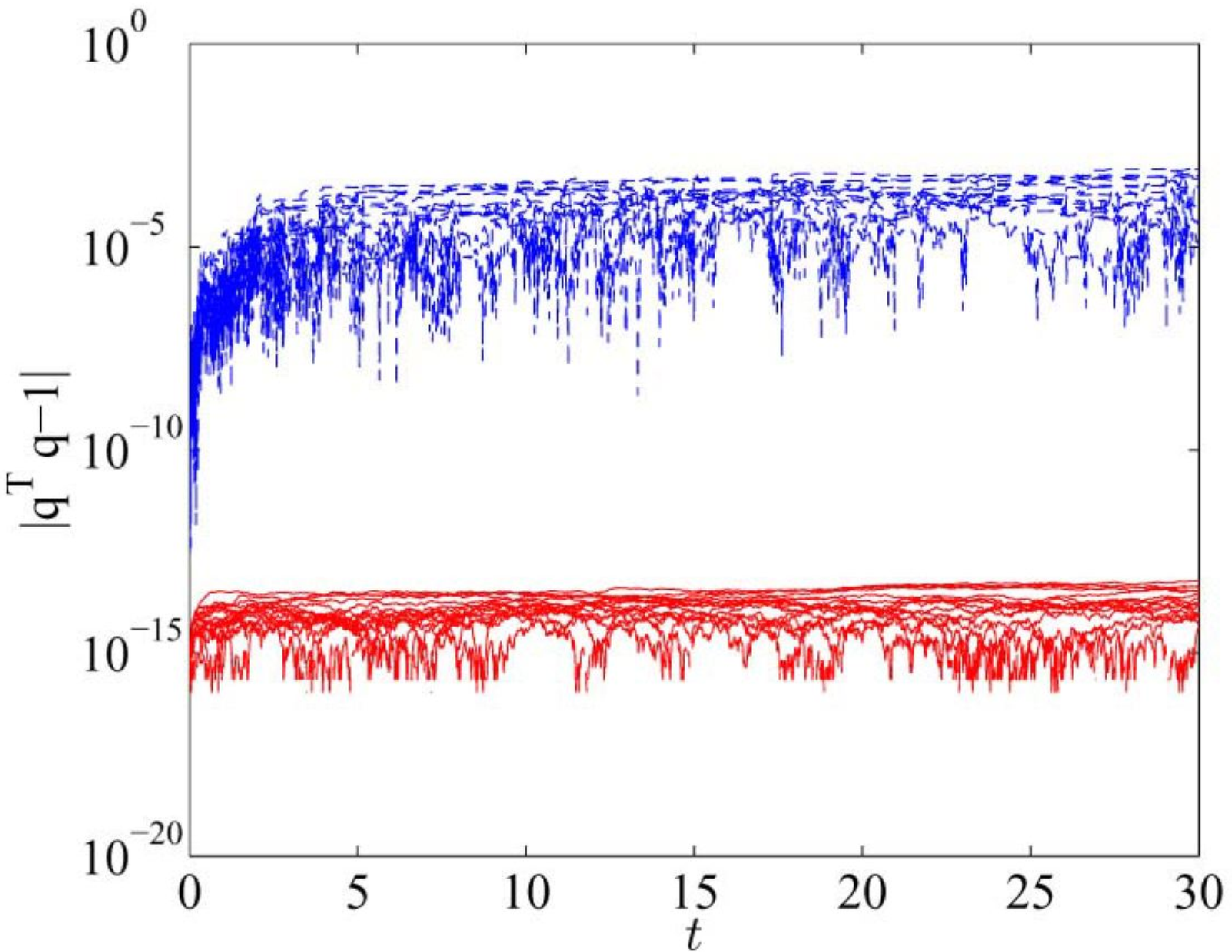}}}
\caption{Numerical simulation of an array of magnetic dipoles (RK45: blue, dotted, VI: red, solid)}\label{fig:mag}
\end{figure}

\reffig{mag} shows the computed total energy and the unit length errors. The variational integrator preserves the total energy and the structure of $(\Sph^2)^n$ well. The mean total energy variation is $8.5403\times 10^{-10}\,\mathrm{Nm}$, and the mean unit length error is $1.6140\times 10^{-14}$. There is a notable dissipation of the computed total energy for the Runge-Kutta method, where the mean variation of the total energy is $2.9989\times 10^{-7}\,\mathrm{Nm}$. The Runge-Kutta method also fail to preserve the structure of $(\Sph^2)^n$. The mean unit length error is $1.7594\times 10^{-4}$.
\end{example}

\begin{example}[\textbf{Molecular Dynamics on a Sphere}]
We study molecular dynamics on $\Sph^2$. Each molecule is modeled as a particle moving on  $\Sph^2$. Molecules are subject to two distinct forces: an attractive force at long range and a repulsive force at short range. Let $m_i\in\Re$ and $q_i\in\Sph^2$ be the mass and the position vector of the $i$-th molecule, respectively. The $i,j$-th element of the inertia matrix is $M_{ij}=m_i$ when $i=j$, and $M_{ij}=0$ otherwise. The Lennard-Jones potential is a simple mathematical model that represents the behavior of molecules~\cite{Len.PPS31}
\begin{align*}
    V(q_1,\ldots,q_n) = \frac{1}{2}\sum_{\substack{i,j=1\\j\neq i}}^n 4\epsilon \bracket{\parenth{\frac{\sigma}{\norm{q_i-q_j}}}^{12}-\parenth{\frac{\sigma}{\norm{q_i-q_j}}}^{6}},
\end{align*}
where the first term models repulsion between molecules at short distance according to the Pauli principle, and the second term models attraction at long distance generated by van der Walls forces. The constant $\epsilon$ and $\sigma$ are molecular constants; $\epsilon$ is proportional to the strength of the mutual potential, and $\sigma$ characterize inter-molecular force. Substituting these into \refeqn{EL}, the continuous equations of motion for the molecular dynamics on a sphere are given by
\begin{align}
    m_i \ddot q_i = -m_i (\dot q_i\cdot \dot q_i)q_i - q_i\times\big(q_i\times\sum_{\substack{j=1\\j\neq i}}^n
    4\epsilon \frac{q_i-q_j}{\norm{q_i-q_j}}\bracket{\frac{12 \sigma^{12}}{\norm{q_i-q_j}^{13}}-\frac{6\sigma^6}{\norm{q_i-q_j}^7}}
    \big)\label{eqn:ELmd}
\end{align}
for $i\in\{1,\ldots,n\}$.

\begin{figure}[t]
\centerline{
    \subfigure[Initial trajectories]{%
        \includegraphics[width=0.22\textwidth]{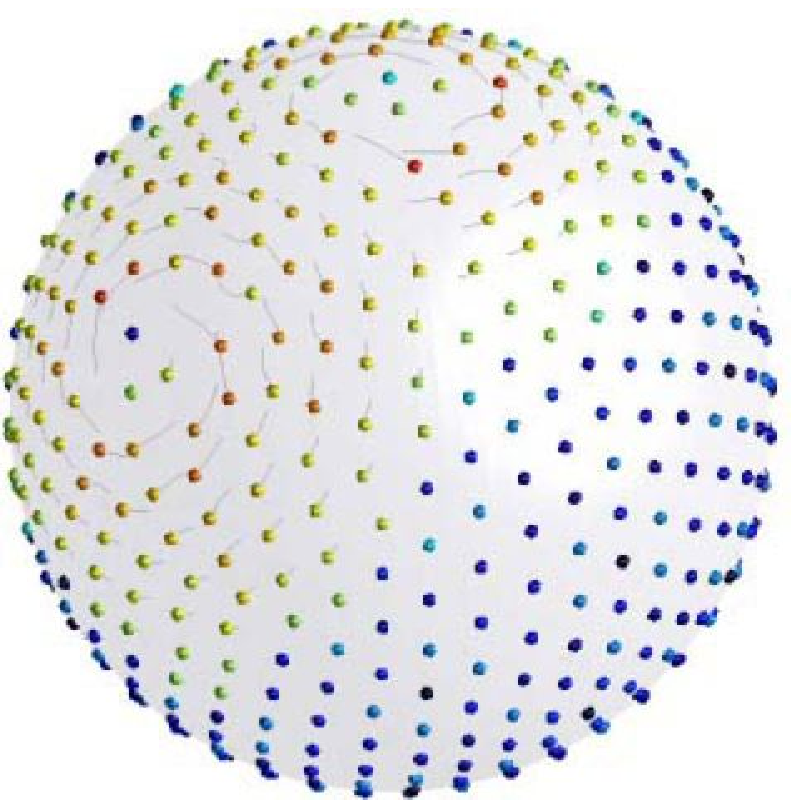}\label{fig:md3d}}
    \hspace{0.02\textwidth}
    \subfigure[Computed total energy]{%
        \includegraphics[width=0.28\textwidth]{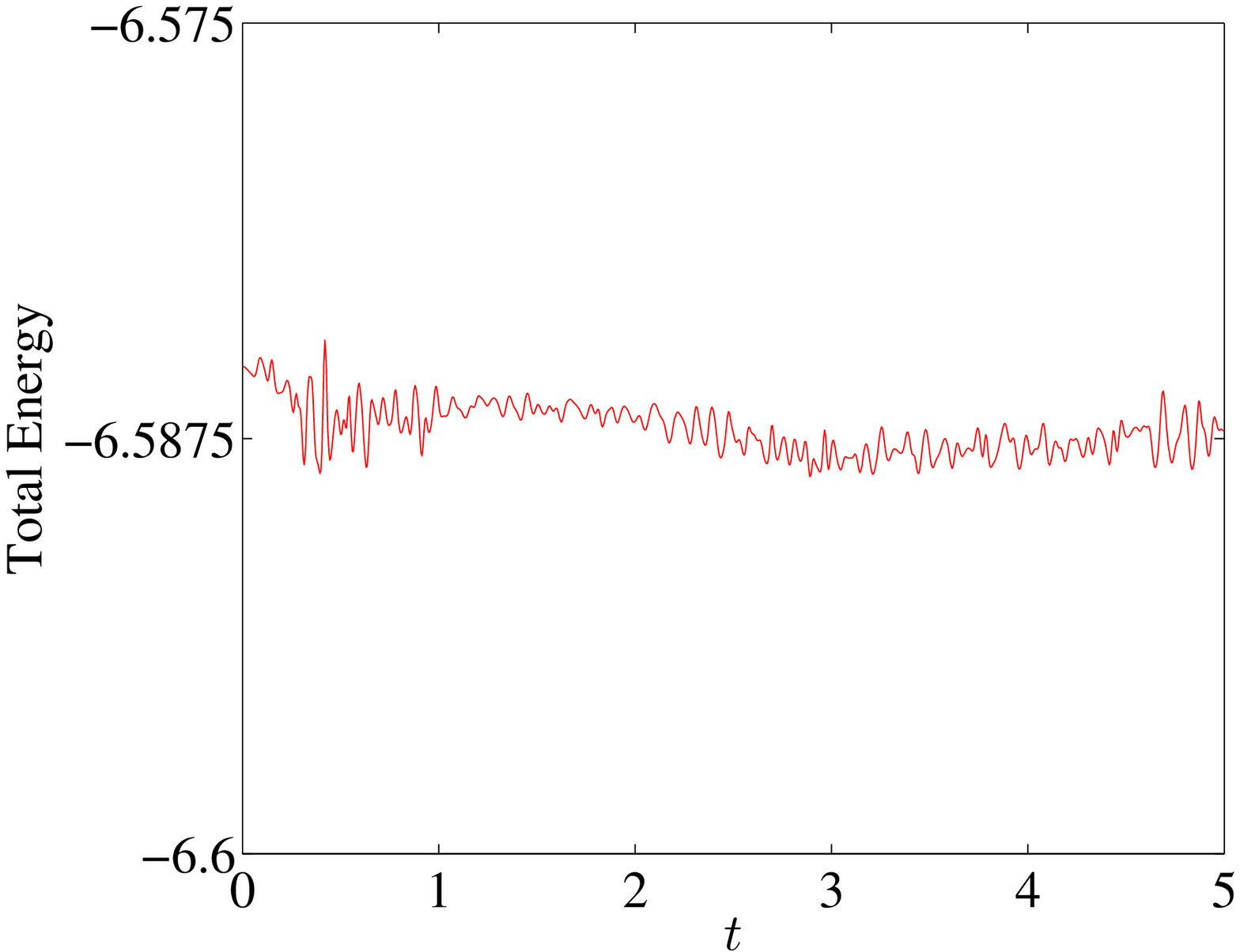}}}
\caption{Numerical simulation of molecular dynamics on a sphere}\label{fig:md}
\end{figure}
\begin{figure}[t]
\centerline{
    \subfigure[$t=0$]{%
        \includegraphics[width=0.195\textwidth]{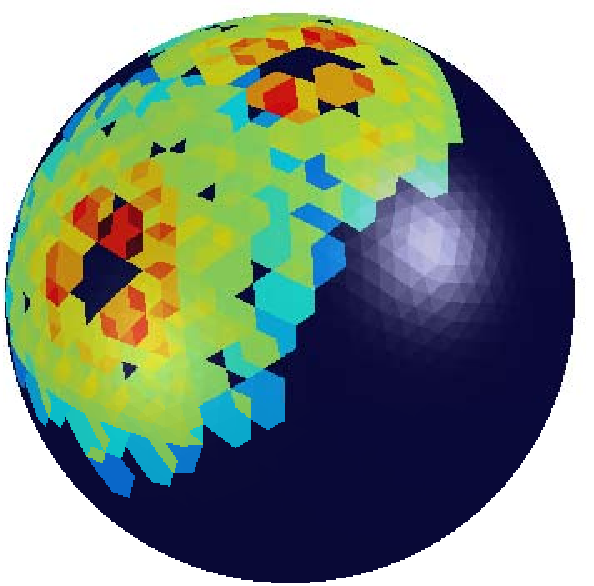}}
    \subfigure[$t=0.25$]{%
        \includegraphics[width=0.195\textwidth]{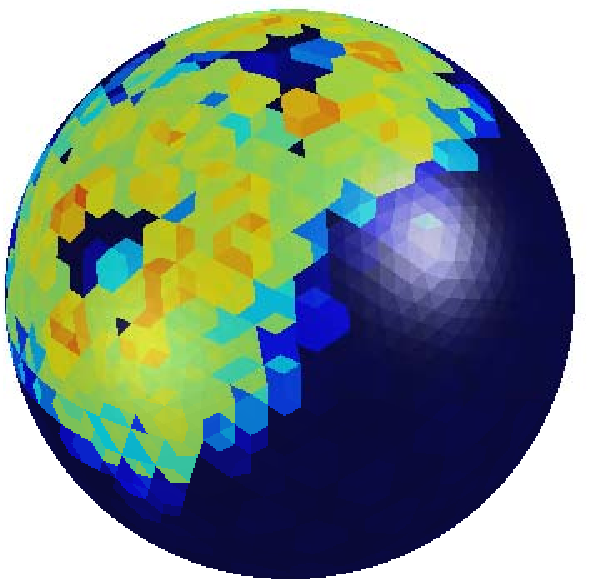}}
    \subfigure[$t=0.5$]{%
        \includegraphics[width=0.195\textwidth]{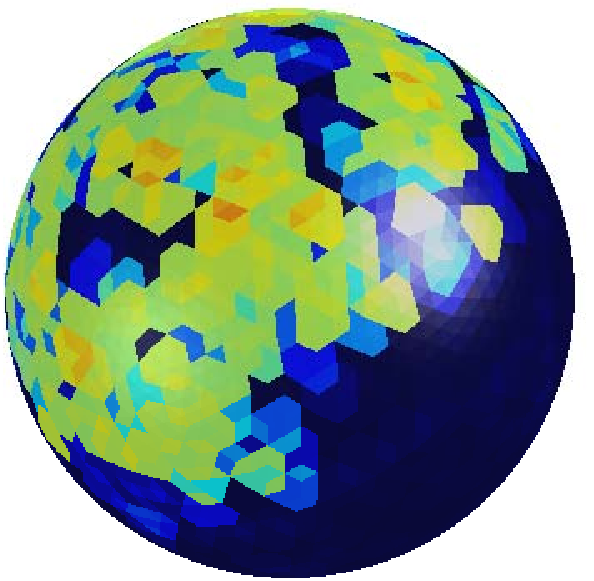}}
    \subfigure[$t=0.75$]{%
        \includegraphics[width=0.195\textwidth]{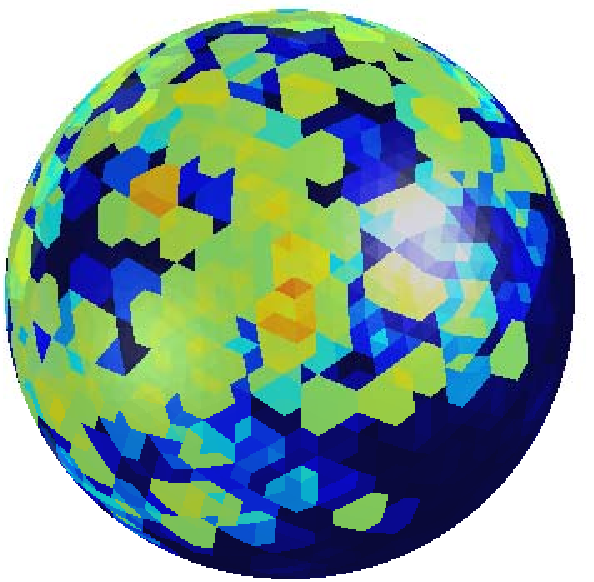}}
    \subfigure[$t=5$]{%
        \includegraphics[width=0.195\textwidth]{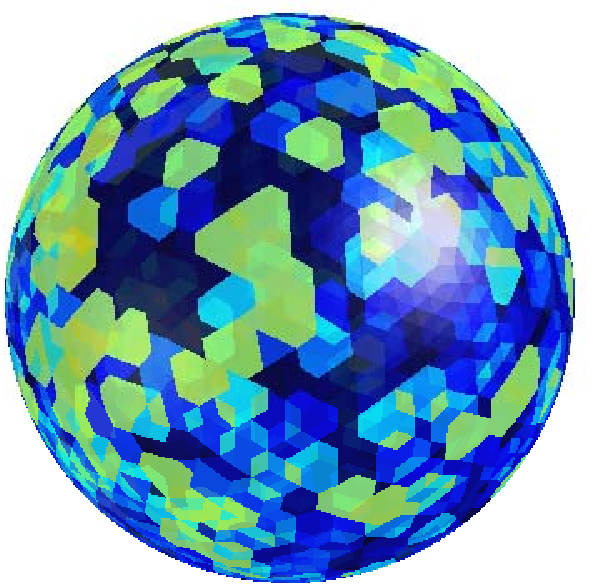}}}
\caption{Kinetic energy distributions over time}\label{fig:mdtemp}
\end{figure}

We choose $642$ molecules, $n=642$, and we assume each molecule has the same mass, $m_i=1$. Initially, molecules are uniformly distributed on a sphere. The strength of the potential is chosen as $\epsilon=0.01$, and the constant $\sigma$ is chosen such that the inter-molecular force between neighboring molecules is close to zero. The initial velocities are modeled as two vortices separated by $30\,\mathrm{degrees}$. The simulation time is $5\,\mathrm{sec}$, and the step size is $h=0.005$.

Trajectories of molecules and the computed total energy is shown at \reffig{md}. The mean deviation of the total energy is $1.8893\times 10^{-3}$, and the mean unit length error is $5.2623\times 10^{-15}$. In molecular dynamics simulations, macroscopic quantities such as temperature and pressure are more useful than trajectories of molecules. \reffig{mdtemp} shows the change of kinetic energy distributions over time, which measures the temperature~\cite{AllTil.BK87}; the sphere is discretized by an icosahedron with $5120$ triangular faces, and the color of a face is determined by the average kinetic energy for molecules that lie within the face and within its neighboring faces. The local kinetic energy is represented by color shading of blue, green, yellow, and red colors in ascending order.

\end{example}

\section{Conclusions}
Euler-Lagrange equations and variational integrators are developed for Lagrangian mechanical systems evolving on $(\Sph^2)^n$ where the Lagrangian is written in a particular form given by \refeqn{L}. The structure of $\Sph^2$ is carefully considered to obtain global equations of motion on $(\Sph^2)^n$ without local parameterizations or explicit constraints.

In the continuous time setting, this provides a remarkably compact form of the equations of motion compared to the popular angular description. For example, it is not practical to study a triple spherical pendulum by using angles due to the complexity of the trigonometric expressions involved. The global Euler-Lagrange equations on $(\Sph^2)^n$ maintain the same compact structure for arbitrary $n$. In particular, it is possible to use them as a finite element model for a continuum problem as shown in Example \ref{ex:rod}. They are also useful for the theoretical study of global dynamic characteristics.

The variational integrators on $(\Sph^2)^n$ preserve the geometric properties of the dynamics as well as the structure of the configuration manifold concurrently. They are symplectic, momentum preserving, and they exhibit good energy behavior for exponentially long time as they are derived from discrete Hamilton's principle. Using the characteristics of the homogeneous manifold $(\Sph^2)^n$, the discrete update map is represented by a group action of $\SO$ to obtain compactly represented, and possibly explicit, numerical integrators. In particular, variational integrators on $(\Sph^2)^n$ completely avoid the singularities and complexity introduced by local parameterizations and explicit constraints.

\bibliographystyle{siam}
\bibliography{gvihs}

\end{document}